\documentclass[11pt]{amsart}
\usepackage[a4paper]{anysize}\marginsize{2.5cm}{2.5cm}{1.3cm}{2cm}
\pdfpagewidth=\paperwidth \pdfpageheight=\paperheight
\usepackage{amsfonts,amssymb,amsthm,amsmath,eucal}
\usepackage{etex}
\usepackage{bbm,array}
\usepackage{graphicx}
\usepackage{pgf,tikz,pgfplots}
\pgfplotsset{compat=1.14}
\usepackage{mathrsfs}
\usepackage[bb=boondox]{mathalfa}
\usetikzlibrary{arrows}
\usetikzlibrary[patterns]
\usepackage{float}
\restylefloat{table}
\usepackage{subfigure}
\usepackage{caption}
\usepackage{setspace}
\usetikzlibrary{arrows}
\usepackage{enumitem}
\usepackage{wrapfig}

\setlist[itemize,3]{leftmargin=-.5in}
\usepackage{titlesec}
\usepackage[normalem]{ulem}
\titlelabel{\thetitle.\;\;}
\usepackage{datetime}
\newdateformat{monthyeardate}{%
	\monthname[\THEMONTH], \THEYEAR}

\theoremstyle{plain}
\newtheorem{thm}{Theorem}[section]
\newtheorem{theorem}[thm]{Theorem}

\newtheorem*{theoremA}{Theorem A}
\newtheorem*{theoremB}{Theorem B}

\newtheorem{lemma}[thm]{Lemma}

\newtheorem{corollary}[thm]{Corollary}

\newcommand\beginproof[1]{\trivlist\item[\hskip\labelsep{\em #1.}]}
\renewcommand\proof{\beginproof{Proof}}

\def\endproof{\hspace*{\fill}\endproofsymbol\endtrivlist}

\def\endproofsymbol{\frame{\rule[0pt]{0pt}{6pt}\rule[0pt]{6pt}{0pt}}}

\theoremstyle{definition}
\newtheorem{definition}[thm]{Definition}
\newtheorem{remark}[thm]{Remark}
\newtheorem{example}[thm]{Example}

\newtheorem{thevarthm}[thm]{\varthmname}

\newenvironment{varthm*}[1]{\trivlist\item[]{\bf #1.}\it}{\endtrivlist}

\renewcommand\geq{\geqslant}

\renewcommand\leq{\leqslant}

\newcommand\be{\begin{eqnarray*}}
	\newcommand\ee{\end{eqnarray*}}

\newcommand\R{\mathbb{ R}}

\newcommand\Z{\mathbb{ Z}}

\newcommand\T{\mathbb{ T}}
\newcommand\K{\mathbb{ K}}
\renewcommand\P{\mathbb{ P}}

\newcommand\newop[2]{\def#1{\mathop{\rm #2}\nolimits}}
\newop\ord{ord}
\newop\mult{mult}
\newop\codim{codim}
\newop\Ass{Ass}
\newop\reg{reg}
\newop\areg{areg}
\newop\satdeg{satdeg}
\newop\supp{supp}
\newop\gin{gin}
\newop\ini{in}
\newop\vol{vol}
\newop\sat{sat}
\newop\length{length}
\newop\depth{depth}
\newop\ahp{aHP}
\newop\HP{HP}
\newop\ri{ri}
\newop\supp{supp}
\newop\ahf{aHF}
\newop\HF{HF}
\newop\dHF {\Delta HF}
\newop\conv{conv}
\newop\lcm{lcm}

\newcommand\ibul{I_{\bullet}}
\newcommand\Ibul{I_{\bullet}}
\newcommand\jbul{J_{\bullet}}

\def\keywordname{{\bfseries Keywords}}%
\def\keywords#1{\par\addvspace\medskipamount{\rightskip=0pt plus1cm
		\def\and{\ifhmode\unskip\nobreak\fi\ $\cdot$
		}\noindent\keywordname\enspace\ignorespaces#1\par}}
\def\subclassname{{\bfseries Mathematics Subject Classification
		(2000)}\enspace}
\def\subclass#1{\par\addvspace\medskipamount{\rightskip=0pt plus1cm
		\def\and{\ifhmode\unskip\nobreak\fi\ $\cdot$
		}\noindent\subclassname\ignorespaces#1\par}}

\definecolor{cqcqcq}{rgb}{0.75294,0.752941,0.75294}

\begin{document}
		
		\author{Grzegorz Malara}
		\title{Asymptotic invariants of generic initial ideals}
		\date{\today}
		\maketitle
		\thispagestyle{empty}
		
		\begin{abstract}
		Generic initial ideals (gins in short) were systematically introduced by Galligo in 1974 under the name of Grauert invariants since they appeared apparently first in works of Grauert and Hironaka. Ever since they are of interest in commutative algebra and indirectly in algebraic geometry. Recently Mayes in a series of articles associated to gins geometric objects called limiting shapes. The construction resembles that of Okunkov bodies but there are some differences as well. This work is motivated by Mayes articles and explores the connections between gins, limiting shapes and some asymptotic invariants of homogeneous ideals, e.g. asymptotic regularity, Waldschmidt constant and some new invariants, which seem relevant from geometric point of view.
		
		In this note we generalize Mayes ideas to graded families of ideals. We work out, sometimes surprising, properties of defined objects and invariants. For example we establish the existence of limiting shapes in dimension $2$ which are polygons with an arbitrary high number of edges and vertices with irrational coordinates.
			
			\vspace{0.2cm}
			\subclass{13P10, 14N20}
		\end{abstract}

\setstretch{1.2}
\section{Introduction}
\thispagestyle{empty} 

In $1974$ Galligo in his paper \cite{Gal} observed that given an ideal $I\subseteq S(n)=\K[x_0,\ldots,x_n]$, general linear change of coordinates $g$ leads to a monomial ideal $\ini(g(I))$ with many interesting properties. 
More specifically, Galligo proved that in $\text{GL}_{n+1}(\mathbb{\K})$ there exists a Zariski open subset $U$, such that the initial ideal $\ini(g(I))$ is invariant for all $g \in U$.
This ideal is thus well-defined, and is called the \emph{generic initial ideal of $I$}, denoted in this paper by $\gin(I)$. The first name that Galligo was using for that kind of ideals was \emph{Grauert invariant}, but since $\gin$'s of ideals are closed under the action of Borel group, a few years later the name of \emph{Borel-fixed ideals} was in use. In characteristic zero the property of being Borel-fixed can be expressed as a nice divisibility property (see Lemma \ref{th:borel_equivalent}) which justifies using the name of \emph{strongly stable ideals}. Fundamental properties of these ideals were established about $20$ years ago in \cite{Green1} and \cite{Green2}.

Mayes \cite{Mayes} introduced geometric figures, limiting shapes $\Delta(\ibul)$ associated to the family $\ibul=(I^{(m)})$  of symbolic powers of a fixed ideal $I$. Dumnicki, Szpond and Tutaj-Gasi\'nska \cite{DST15} for symbolic and ordinary powers of ideals introduced the notions of asymptotic Hilbert function and asymptotic Hilbert polynomial. Ordinary and symbolic powers of an ideal form a graded sequence of ideals. Thus it is natural to try to extend the ideas of limiting shapes and asymptotic Hilbert invariants to arbitrary graded sequences of ideals. Relevant constructions and their properties fill the core of the paper. We provide also some non-trivial applications in Section \ref{sec:limiting shape} and \ref{sec:0dim}. 

It was observed by Macaulay that the initial ideal $\ini(I)$ of an ideal $I$ encodes almost all information on $I$. The advantage of working with the generic initial ideal $\gin(I)$ is that, in contrary to $\ini(I)$, it does not depend on the particular choose of coordinates. In the 90's first asymptotic ideas of ideals (regularity of ordinary and symbolic powers of an ideal) appear in works of Geramita, Gimigliano, Pitteloud \cite{GGP} and Chandler \cite{Chandler}. Combining the asymptotic approach with gins and associated geometrical objects (limiting shapes) reveals certain additional invariants. For example, it is known that for graded sequence of symbolic or ordinary powers of ideal $\ibul$ the Waldschmidt constant $\widehat{\alpha}(\ibul)$ is the minimal first coordinate of intersection points of the limiting shape and the first coordinate axis. The same happens for any graded sequence of ideals, namely there is
$$\widehat{\alpha}(\ibul)=\min\Big\{\alpha_0\;:\;(\alpha_0,\alpha_1,\ldots,\alpha_n)\in \Big(\Delta(\ibul,t)\cap \{x_1=\ldots=x_n=0\}\Big)\subseteq \R^{n+1}\Big\}. $$ For the asymptotic regularity $\widehat{\reg}(\ibul)$ the situation is more involved. We prove that
$$\widehat{\reg}(\ibul)=\sup \{|x|, \text{ such that } x=(x_0,x_1,\ldots,x_n) \text{ is an extremal point in }\Delta(\ibul)\} \in \R\cup \{\infty\},$$
where $|x|=|x_0|+\ldots+|x_n|$. In particular $\widehat{\reg}(\ibul)$ may not be finite (see Example \ref{ex:ahp=1}). We provide examples of graded sequences of ideals with irrational Waldschmidt constant and irrational asymptotic regularity, see Example \ref{ex:Waldschmidt areg irrational}. It remains an open problem if these invariants may be irrational for symbolic powers of an ideal.

The main result of this paper is the construction of limiting shapes with arbitrary high number of line segments which form its boundary. It is done for symbolic powers of ideals defined by carefully chosen $0$-dimensional schemes in $\P^2$ (see Corollary \ref{cor:n break points gin}). 
\begin{theoremA}
	For any positive integer $M>0$, there exists a set of points $Z$ in $\P^2$ such that the limiting shape of the symbolic powers of ideal $I(Z)$ has at least $M$ line segments.
\end{theoremA}
We prove also the following related result (see Example \ref{ex:anyDescent})
\begin{theoremB}
	For a positive integer $M>0$, there exists a graded sequence of ideals $\ibul$ (depending on $M$) such that its limiting shape has at least $M$ vertices and all coordinates of these vertices are irrational numbers.
\end{theoremB}
It is not clear if the vertices of limiting shapes can be explained by geometry of the schemes defined by ideals in the sequence. In the specific Corollary \ref{cor:n break points gin} justifying Theorem A, we expect that the vertices are related to collinearity of points in $Z$ but we were not able to make this relation precise.

This paper is organized as follows. For a graded sequence of ideals we define asymptotic Hilbert function and asymptotic Hilbert polynomial in Section \ref{sec:functions} using sets constructed in Section \ref{sec:various sets}. Definitions of limiting shape and its complement are formulated in Section \ref{sec:limiting shape}. This section ends with Example \ref{ex:anyDescent}, which justifies Theorem B.  In the next section we introduce the first difference asymptotic Hilbert function, which we use in the last Section \ref{sec:0dim}, where we turn our attention to $0$-dimensional schemes in $\P^2$ and families of symbolic powers of their ideals. In this Section \ref{sec:0dim} we show the applications of developed theory, culminating in Theorem \ref{th:n lines dHF} and Theorem \ref{th:2lines+1pts dHF} from which Theorem A follows.

\section{Various sets associated to homogeneous ideals}
\label{sec:various sets}
We begin by studying the action general linear group on the ring of polynomials $S(n)=\K[x_0,\ldots,x_n]$ over a field $\K$ of characteristic $0$.

Consider the general linear group $\text{GL}_n(\K)$ as the set of invertible matrices $g=(c_{ij}) \in \text{GL}_n(\K)$. Then for a polynomial $f \in S(n-1)$ we have
$$g(f(x)):=f(g(x))=f\Bigg(\sum_{j=0}^{n-1} c_{1j}x_j,\ldots,\sum_{j=0}^{n-1} c_{nj}x_j\Bigg). $$
In general, for an ideal $I \subset S(n-1)$, we define $ g(I)=\{g(f)\; : \; f \in I\}.$
We say that an ideal $I$ is \textit{fixed under} a subgroup $G \subset \text{GL}_n(\K)$ (or simply $G$--\textit{fixed}) if $g(I) \subset I$ for all $g \in G$.

In $\text{GL}_n(\K)$ we distinguish the subgroup of upper triangular matrices, denoted by $\text{B}_n(\K)$, the so-called Borel group, for which we have the following nice property.
\begin{lemma}
	\label{th:borel_equivalent}
	For a monomial ideal $I \subset S(n-1)$ the following conditions are equivalent
	\begin{itemize}
		\item [1)] $I$ is $\text{B}_n$-fixed;
		\item [2)] If $\omega \in I$ is any monomial divisible by $x_j$, then $\omega\frac{x_i}{x_j} \in I$ for all $i<j$.
	\end{itemize}
\end{lemma}

In this paper we consider a graded sequence of ideals, i.e. a family of ideals $\ibul=\{I_m\}_{m \in \Z_{\geq 0}}$ such that for all $p,q \in \Z_{\geq 0}$ we have $I_p \cdot I_q \subseteq I_{p+q}$. Easy examples of graded sequences of ideals are given by ordinary powers or symbolic powers of a fixed ideal.

Denote by $\ini(f)$ the initial term of $f$ in fixed but arbitrary term order, and by $\ini(I)$ the initial ideal of $I$. Let $\ibul$ be a graded sequence of ideals. The ideal $\ibul^{\text{in}}:=\ini(I_k)_{k\geq 1}$ form also a graded sequence of ideals.

Using the identification $\ini(f)=x^{\alpha}\rightarrow \alpha=(\alpha_0,\ldots,\alpha_n)$ to every polynomial from $f \in S(n)$ we can assign a point in $\R^{n+1}$. For a real number $s$ we consider sets
$$\T_s^n:= \{\beta \in \R^n_{\geq0} :  \; \beta_0+ \ldots +\beta_{n-1} \leq s \}.$$

From now on $m$ will always be a non-negative integer and $t$ a real number.

For a graded sequence of homogeneous ideals $\ibul$ in  $S(n)$ we define sets
$$M_{m,t}(\ibul^{\ini}):=\bigcup_{x^{\alpha} \in \ini(I_m)}\{\beta \in \R^n : \beta_j\geq \alpha_j\; \text{ for all } j=0,1,\ldots,n-1\text{ and } \;\beta_0+ \ldots +\beta_{n-1} \leq mt -\alpha_{n} \}$$
and
$$L_{m,t}(\ibul):=\Big(\bigcup_{\alpha \in M_{m,t}(\ibul^{\ini })}\{\alpha + \R^n_{\geq0} \}\Big) \cap \T_{mt}^n. $$

Note that the sets $L_{m,t}(\ibul)$ fulfil the condition
\begin{equation}
\label{eq:Lmt minkowski}
	L_{p,t}(\ibul) + L_{q,t} (\ibul)\subset L_{p+q,t}(\ibul), 
\end{equation}
where the sum on the left is the Minkowski sum of sets.

Indeed, it follows from the fact that $\ibul$ is a graded sequence of ideals. For the sets $L_{m,t}(\ibul)$ we can prove even more.

\begin{lemma}
	\label{th:Lmt HomogeneousConvex}
	For a graded sequence of homogeneous ideals $\Ibul$ the set $ U:=\overline{\bigcup \frac{L_{m,t}(\ibul)}{m}}$ is convex.
\end{lemma}
\proof
	Take $\alpha,\beta \in U$. There are $\{\alpha_i\}_i,\{\beta_j\}_j$ such that $\alpha_i \in  \frac{L_{p_i,t}(\ibul)}{p_i}$ and $\beta_j \in  \frac{L_{q_j,t}(\ibul)}{q_j}$, and $\lim _{i \rightarrow \infty}  \alpha_i =\alpha$, $\lim _{j \rightarrow \infty}  \beta_j =\beta$.
	
	Let $\gamma=\lambda\alpha+ (1-\lambda) \beta $ for $\lambda \in [0,1 ]$. Let $\lambda_i =\frac{k_i}{s_i}$ with $ k_i,s_i \in \Z_{\geq 0}$ be rational numbers with $\lim _{i \rightarrow \infty}  \lambda_i =\lambda$.
	Then we have
	$$\gamma_i=\lambda_i\alpha_i+ (1-\lambda_i) \beta_i = \frac{1}{s_i}\Big(k_i \alpha_i + (s_i - k_i) \beta_i \Big) \in \frac{1}{s_i} \Big( \frac{L_{k_i p_i,t}(\ibul)}{p_i} + \frac{L_{(s_i-k_i)q_i,t}(\ibul)}{q_i} \Big) \subset$$
	$$\subset \frac{1}{s_i} \frac{L_{k_i p_i q_i+(s_i-k_i)q_i p_i,t}(\ibul)}{p_i q_i} = \frac{1}{s_i} \frac{L_{s_i p_i q_i,t}(\ibul)}{p_i q_i} \subset  \bigcup_{m \geq 1} \frac{L_{m,t}(\ibul)}{m}.$$
	Passing with $i$ to the infinity we obtain $\gamma \in U$.
\endproof

\begin{theorem}
	\label{th:Lmt_sequence}
	For an arbitrary graded family of ideals $\ibul$ and a fixed real number $t$ the limit
	$$\lim_{m \rightarrow \infty} \vol\Big( \frac{L_{m,t}(\ibul)}{m}\Big) = \sup_m \vol \Big( \frac{L_{m,t}(\ibul)}{m}\Big)=\vol \Big(\bigcup_m \frac{L_{m,t}(\ibul)}{m}\Big)$$
	exists.
\end{theorem}
\proof
	The proof of this theorem is similar to proof of \cite[Theorem 5]{DST15}, since it depends only on property (\ref{eq:Lmt minkowski}).
\endproof

\begin{definition}
	\label{def:gamma_Homogeneous}
	Let $t$ be a given real number and $\Ibul$ a graded sequence of homogeneous ideals. For sets $L_{m,t}(\ibul)$ we define
	$$ \Gamma _{m,t}(\ibul) := \T_{mt}^n \setminus L_{m,t}(\ibul).$$
\end{definition}
For a set $\Omega \subseteq \R^n$ we denote by $\# \Omega$ the number of points in $\Omega \cap \Z^n$.
\begin{theorem}
	\label{th:gamma_vol homogeneous}
	For a graded sequence of homogeneous ideals such that for each $m$ $\ini(I_m)$ is a Borel-fixed ideal, the following limits exist and are equal $$ \lim _{m \rightarrow \infty}  \frac{\# \Gamma_{m,t}(\ibul)}{m^{n}} = \lim _{m \rightarrow \infty} \vol \Big( \frac{ \Gamma_{m,t}(\ibul)}{m} \Big).$$
\end{theorem}
\proof
	Observe that the equality $$ \lim _{m \rightarrow \infty} \vol \Big( \frac{ \Gamma_{m,t}(\ibul)}{m} \Big)  = \frac{t^n}{n!}- \lim _{m \rightarrow \infty} \vol \Big( \frac{ L_{m,t}(\ibul)}{m} \Big) $$
	is a consequence of Definition \ref{def:gamma_Homogeneous}.
	Define the sets
	$$A_m:= \Big\{ (c_0,\ldots,c_{n-1}) \in \Z^{n}_{\geq0} : (c_0+1,\ldots,c_{n-1}+1) \in \overline{\Gamma_{m,t}(\ibul)}\Big\}, $$
	$$B_m:=(\Gamma_{m,t}(\ibul) \cap \Z^{n}) \setminus A_m , $$
	which in fact means that $\# \Gamma_{m,t} = \#A_m+ \# B_m$.
	
	Observe that with a point $a=(a_0,\ldots,a_{n-1}) \in A_m$, the cube 
	$$a + C^{[n]}:=\underbrace{[a_0,a_0+1] \times \ldots \times [a_{n-1},a_{n-1}+1]}_{\text{$n$}}  \subset \R^n$$ is fully contained in $\Gamma_{m,t}(\ibul)$. Therefore
	$$\vol\Big(\Gamma_{m,t}(\ibul)\Big) = \# A_m + \vol\Big( (B_m)_C \Big), $$
	where $(B_m)_C= \Gamma_{m,t}(\ibul) - \bigcup_{a\in A_m} \Big( a+ C^{[n]} \Big)$. Note that $(B_m)_C \cap \Z^n=B_m$.
	
	\begin{minipage}{0.45\textwidth}
	Before we proceeding with the proof, it might be convenient to consider the following example and illustration of introduced sets. Let $I=(xy^6,x^3y^5z,x^2yz^3,x^4z)$, then $L_{1,9}(I)$ is the shaded area indicated in the presented figure. Here the set $T_9^2$ is the triangle $(0,0)$, $(0,9)$ and $(9,0)$. Points in the sets $A$ are denoted by filled circle, while points from the set $B$  are denoted as crosses.
	\end{minipage}
		\begin{minipage}{0.45\textwidth}
			\begin{figure}[H]
				\label{fig:AB L19}
				\centering
				\begin{tikzpicture}[line cap=round,line join=round,>=triangle 45,x=1.0cm,y=1.0cm,scale=0.45]
				\draw [color=cqcqcq,, xstep=1.0cm,ystep=1.0cm] (-0.48,-0.5) grid (9.42,9.46);
				\draw[->,color=black] (0.,0.) -- (9.42,0.);
				\foreach \x in {,2.,4.,6.,8.}
				\draw[shift={(\x,0)},color=black] (0pt,2pt) -- (0pt,-2pt) node[below] {\footnotesize $\x$};
				\draw[->,color=black] (0.,0.) -- (0.,10.46);
				\foreach \y in {,2.,4.,6.,8.}
				\draw[shift={(0,\y)},color=black] (2pt,0pt) -- (-2pt,0pt) node[left] {\footnotesize $\y$};
				\draw[color=black] (0pt,-10pt) node[right] {};
				\clip(-0.48,-0.5) rectangle (9.42,9.46);
				\fill[fill=black,pattern=north east lines,pattern color=black] (1.,8.) -- (1.,6.) -- (3.,6.) -- cycle;
				\fill[fill=black,pattern=north east lines,pattern color=black] (2.,7.) -- (2.,1.) -- (8.,1.) -- cycle;
				\fill[fill=black,pattern=north east lines,pattern color=black] (4.,5.) -- (4.,0.) -- (9.,0.) -- cycle;
				\fill[fill=black,pattern=north east lines,pattern color=black] (3.,6.) -- (3.,5.) -- (4.,5.) -- cycle;
				\draw [dash pattern=on 2pt off 2pt] (0.,9.)-- (9.,0.);
				\draw (1.,8.)-- (1.,6.);
				\draw (1.,6.)-- (3.,6.);
				\draw (3.,6.)-- (1.,8.);
				\draw (2.,7.)-- (2.,1.);
				\draw (2.,1.)-- (8.,1.);
				\draw (8.,1.)-- (2.,7.);
				\draw (4.,5.)-- (4.,0.);
				\draw (4.,0.)-- (9.,0.);
				\draw (9.,0.)-- (4.,5.);
				\draw (3.,6.)-- (3.,5.);
				\draw (3.,5.)-- (4.,5.);
				\draw (4.,5.)-- (3.,6.);
				\begin{scriptsize}
				\draw [fill=black] (0.,7.) circle (2.5pt);
				\draw [fill=black] (0.,6.) circle (2.5pt);
				\draw [fill=black] (0.,5.) circle (2.5pt);
				\draw [fill=black] (1.,5.) circle (2.5pt);
				\draw [fill=black] (1.,4.) circle (2.5pt);
				\draw [fill=black] (0.,4.) circle (2.5pt);
				\draw [fill=black] (0.,3.) circle (2.5pt);
				\draw [fill=black] (1.,3.) circle (2.5pt);
				\draw [fill=black] (0.,2.) circle (2.5pt);
				\draw [fill=black] (1.,2.) circle (2.5pt);
				\draw [fill=black] (0.,1.) circle (2.5pt);
				\draw [fill=black] (1.,1.) circle (2.5pt);
				\draw [fill=black] (0.,0.) circle (2.5pt);
				\draw [fill=black] (1.,0.) circle (2.5pt);
				\draw [fill=black] (2.,0.) circle (2.5pt);
				\draw [fill=black] (3.,0.) circle (2.5pt);
				\draw [color=black] (0.,8.)-- ++(-3pt,-3pt) -- ++(6.0pt,6.0pt) ++(-6.0pt,0) -- ++(6.0pt,-6.0pt);
				\draw [color=black] (0.,9.)-- ++(-3pt,-3pt) -- ++(6.0pt,6.0pt) ++(-6.0pt,0) -- ++(6.0pt,-6.0pt);
				\end{scriptsize}
				\end{tikzpicture}
		\end{figure}
	\end{minipage}

	We need to show that
	$$\lim _{m \rightarrow \infty}  \frac{\vol(B_m)_C}{m^{n}}=0.$$
	First we claim that for every point $b=(b_0,\ldots,b_{n-1})$ from the set $B_m$ we have
	\begin{equation}
	\label{eq:star}
		\lfloor mt \rfloor-n \leq  \sum_{i=0}^{n-1} b_i \leq \lfloor mt \rfloor.\tag{*}
	\end{equation}
	In order to show it, take any $\gamma=(\gamma_0,\ldots,\gamma_{n-1}) \in \Gamma_{m,t}(\ibul) \cap \Z^n$. From definition of $\Gamma_{m,t}(\ibul)$ and the fact that all $\gamma_i$ are non--negative integers, we have
	$$\gamma_0+\ldots+\gamma_{n-1} \leq \lfloor mt \rfloor .$$
	
	Assume that $\gamma_0+\ldots+\gamma_{n-1} < \lfloor mt \rfloor -n $. We will show that this implies $\gamma \in A_m$.
	Take the point $ \gamma + \mathbb{1}:=(\gamma_0 +1,\ldots,\gamma_{n-1}+1)$. Then $\gamma + \mathbb{1} \in \T_{mt}^n$ since we have
	$(\gamma_0+1)+\ldots+(\gamma_{n-1}+1) \leq \lfloor mt \rfloor$. Assume that $\gamma + \mathbb{1} \in L_{m,t}(\ibul)$, then we have two possibilities.
	\begin{itemize}
		\item [1)] If $\gamma + \mathbb{1} \in \partial(L_{m,t}(\ibul))$, then $\gamma + \mathbb{1}  \in \partial(\Gamma_{m,t}(\ibul))$, hence $\gamma \in A_m$.
		
		\item [2)] If $\gamma + \mathbb{1} \in \text{int}(L_{m,t}(\ibul))$, then  $\gamma \in L_{m,t}(\ibul)$, a contradiction to the assumption $\gamma \in \Gamma_{m,t}(\ibul)$.
	\end{itemize}
This shows the claim (\ref{eq:star}).
	Turning back to the main proof, we calculate the number of integral points $b$ fulfilling the condition
	$$\lfloor mt \rfloor+n \geq b_0+1+\ldots+b_{n-1}+1=b_0+\ldots+b_{n-1}+n\geq \lfloor mt \rfloor.$$
	These points lie in hyperplanes described by equations
	$$x_0+\ldots+x_{n-1}=\lfloor mt \rfloor-n,\;\;x_0+\ldots+x_{n-1}=\lfloor mt \rfloor-n+1,\;\ldots,\;x_0+\ldots+x_{n-1}=\lfloor mt \rfloor .$$
	The number of these points is not larger than the maximal possible number of non--negative integer points satisfying these equations, i.e.,
	$$\binom{\lfloor mt \rfloor-1}{n-1}+\binom{\lfloor mt \rfloor}{n-1}+\ldots+\binom{\lfloor mt \rfloor+n-1}{n-1} \leq (n+1)\binom{\lfloor mt \rfloor+n-1}{n-1}.$$
	
	Thus we have
	$$0\leq \#B_m\leq (n+1)\binom{\lfloor mt \rfloor+n-1}{n-1}=(n+1)\frac{(\lfloor mt \rfloor+1)\cdot\ldots\cdot(\lfloor mt \rfloor+n-1)}{(n-1)!}, $$
	and then
	$$\lim_{m\rightarrow \infty}\frac{\vol((B_m))_C}{m^n}=0. $$
\endproof

\section{Functions associated to Borel-fixed monomial ideals}
\label{sec:functions}

Now we extend classical Hilbert function for ideal $I$ by setting $$\widetilde{\HF}_{I}(t):=\HF_{I}(\lfloor t\rfloor)=\dim_{\K}(S(n) /I)_{\lfloor t\rfloor} . $$
If we use the definition of $\Gamma _{m,t}(\ibul)$ just for one ideal $I$, by putting $m=1$, we get $\widetilde{\HF}_{I}(t):=\# \Gamma_t(I), $ for all $t \in \R$.
\begin{theorem}
	\label{th:ahf_homogeneous_ existence}
	Assume that we have
	\begin{itemize}
		\item [a)] a graded sequence of homogeneous ideals,
		
		or
		
		\item [b)] a graded sequence of ideals with degree compatible term ordering,
	\end{itemize}
	   such that for each $m$ the ideal $\ini(I_m)$ is a Borel-fixed ideal. Then there exists the limit $$\lim_{m \rightarrow \infty} \frac{\widetilde{\HF}_{I_m}(mt)}{m^n}, $$
	for any real number $t$. 
\end{theorem}
\proof
We begin with the proof that $\widetilde{\HF}_{I_m}(mt)=\widetilde{\HF}_{\ini(I_m)}(mt)$ for $t \in \Z$. By Macaulay's Basis Theorem (\cite[Theorem 15.3]{Eis}) the set of all monomials not in $\ini(I_m)$ forms a basis for $S(n)/ I_m$. We need to show that therefore the same holds in degree $t$.

Let $I=(f_1,\ldots,f_k)$ and denote $$B_t=\{ x^a + (I_m)_t\; :\;\;x^a \notin \ini(I_m)_t\; \text{ and }\; |a|=t \}.$$ We prove that $B_t$ is a basis for $(S(n) / I_m)_t$. To show that $B_t$ forms a system of generators let $f \in S(n)$ be such that $\deg(f)=t$. Here our proof splits in two different parts depending on what we assume:

\noindent
a) $I_m$ is a graded sequence of  homogeneous  ideals, then
$$f=f_{(I_m)_t} +(I_m)_t,$$ where $f_{(I_m)_t}$ is the remainder of division by $(I_m)_t$.

\noindent
b) $I_m$ is a graded sequence of ideals with the assumption that the term ordering is degree compatible. This assures that dividing $f$ by $(I_m)_t$ gives us the reminder which does not belong to the quadrant $\{ x_0^{b_1},\ldots,x_n^{b_n}\; : \;  b_i \geq a_i \}$ for any $x^{a_1,\ldots,a_n} \in \ini(I_m)_t$, i.e. $f=f_{(I_m)_t} +(I_m)_t$ as in the previous case.

In both cases, if $\deg(f_{(I_m)_t}) \neq 0$, then $\deg(f_{(I_m)_t})=t$ and $f_{(I_m)_t} \notin \ini(I_m)_t$, which finishes the proof.

To prove the linear independence, suppose that $$\sum_{i=1}^{s} c_i (x^{a_i} +  (I_m)_t ) =0,\;\; x^{a_i} \notin \ini(I_m)_t,\;\;|a_i|=t,\;\; c_i \in \mathbb{K}.$$ But this means that $\sum_{i=1}^{s} c_i x^{a_i} \in (I_m)_t $ and the only possibility is that $$\sum_{i=1}^{s} c_i x^{a_i} =0,\;\; \text{ hence } \;\; c_i=0.$$

The number $\widetilde{\HF}_{\ini(I_m)}(mt)$ is equal to the number of all monomials in $\Big( S(n) /\ini(I_m)\Big)_{mt}$. Therefore by definition of $\Gamma_{m,t}$
 $$\#\Gamma_{m,t} =  \widetilde{\HF}_{\ini(I_m)}(mt).$$

Applying Theorem \ref{th:gamma_vol homogeneous}  completes the proof.
\endproof

Proved theorem allow us to extend the definition of the asymptotic Hilbert function, originally stated in \cite[Definition 1]{DST15}, by the cases of homogeneous graded sequences of ideals.

\begin{definition}
	\label{def:aHF Homogeneous}
	For any real number $t$ and a graded sequence of ideals $\Ibul$ fulfilling the assumptions of Theorem \ref{th:gamma_vol homogeneous}, the \textit{asymptotic Hilbert function} of $\Ibul$ is defined as
	$$\ahf_{\ibul} (t) := \lim_{m \rightarrow \infty} \frac{\widetilde{\HF}_{I_m}(mt)}{m^n}. $$
\end{definition}
To give similar definition for asymptotic Hilbert polynomial $\ahp$ we want to use the regularity index $\ri(I)$ (also called Hilbert regularity), which is the smallest number starting from which $\HF_I(i) =\HP_I(i)$ for all $i\geq \ri(I)$. The existence of $\ahp$ was proven in \cite[Theorem 13]{DST15} for radical homogeneous ideals under the condition that \textit{the linearly bounded symbolic regularity} (see \cite{HerzogHoaTrung}  and \cite{Swanson}) holds. We extend this property to graded sequences of Borel-fixed monomial ideals, using regularity index.

\begin{definition}
	\label{def:LBHR}
	Let $\ibul$ be a graded sequence of monomial ideals in $S(n)$. We say that $\ibul$ satisfies \textit{linearly bounded Hilbert regularity} (or LBHR for short) if there exist constants $a,b >0$, such that $$\ri(I_m) \leq am+b \text{ for all } m \geq 0.$$
\end{definition}
\begin{definition}
	\label{def:aHP}
	The \textit{asymptotic Hilbert polynomial} of a graded sequence of Borel-fixed monomial ideals $\ibul$ is
	$$\ahp_{\ibul}(t) :=\lim_{m \rightarrow \infty} \frac{\HP_{I_m}(mt)}{m^n}. $$
\end{definition}
\begin{theorem}
	\label{th:ahf ahp polynomial monomial}
	For a graded sequence of Borel-fixed monomial ideals $\Ibul$ in $S(n)$, for which linearly bounded Hilbert regularity holds, $\ahp_{\ibul}$ exists, $\ahp_{\ibul}(t)=\ahf_{\ibul}(t)$ for $t \geq a+b$ and $\ahp_{\ibul}$ is indeed a polynomial.
\end{theorem}
\proof
See proof of \cite[Theorem 13 and Theorem 15]{DST15} and take $\ri(I_m)$ instead of $\reg(I_m)$.
\endproof

\section{Asymptotic invariants}
\label{sec:invariants}

We begin this section by adapting well-known asymptotic invariant, so called Waldschmidt constant, to the graded sequence of ideals.
\begin{definition}
	\label{def:WaldschmidtConst}
	Let $0 \neq I\neq S(n)$, be a homogeneous ideal and denote by $\alpha(I):=\min \{d: I_d \neq 0\}$
	the least degree of a non-zero element in $I$.
	The Waldschmidt constant of $\ibul$ is defined as the real number
	$$ \widehat{\alpha}(\ibul) :=\lim_{m \rightarrow \infty} \frac{\alpha(I_m)}{m},$$
	where $I_m \in \ibul$.
\end{definition}
We note here that for a graded system of ideals $\Ibul$ the asymptotic invariant $\widehat{\alpha}(\ibul)$ always exists. We can prove it by using Fekete Lemma \cite{Fekete} since from $I_p \cdot I_q \subset I_{p+q}$ we get $ \alpha(I_{p+q}) \leq \alpha(I_p) + \alpha(I_q) $ for all positive integer $p,q$.

The notation for Waldschmidt constants has evolved over the years. Whereas Waldschmidt, apparently has not used any symbol, Chudnovsky \cite{Chudnovsky} used $\Omega_0(\ibul)$. In commutative algebra Waldschmidt constants have been introduced in works of Bocci, Dumnicki, Harbourne, Szemberg and Tutaj-Gasińska \cite{BoHa1}, \cite{BoHa2}, \cite{DHST}. They used the symbol $\gamma(\ibul)$. The notation $\widehat{\alpha}(\ibul)$, resembling the notation for asymptotic cohomology functions \cite{FKL} has appeared first in \cite{DHST}. It is also worth to mention here that so far, in the literature, Waldschmidt constant is calculated only for symbolic and ordinary powers of ideal, a special cases of graded ideals. 

In the paper \cite{CutkoskyHerzogTrung} Cutkosky, Herzog and Trung introduce asymptotic quantities, the so-called {\it asymptotic regularity} for ordinary powers of ideals. This work had its continuity in the work of Cutkosky and Kurano in \cite{CutkoskyKurano}, where the authors studied this notion for symbolic powers of ideals. Here we introduce the same asymptotic quantities for a graded family of ideals $\ibul$.
\begin{definition}
	\label{def:areg Homogeneous}
	Let $\Ibul$ be a family of homogeneous ideals such that for each $m$ the ideal $\ini(I_m)$ is a Borel-fixed ideal in the polynomial ring $S(n)$, then the asymptotic regularity and the asymptotic Hilbert regularity of $\Ibul$ are define respectively
	$$\widehat{\reg}(\ibul) :=\lim_{m \rightarrow \infty} \frac{\reg(I_m)}{m}, \qquad \widehat{\ri}(\ibul) :=\lim_{m \rightarrow \infty} \frac{\ri(I_m)}{m},$$
	if the limits exist and are finite.
\end{definition}
The existence of $\widehat{\reg}(\ibul)$ was proved for some special cases of graded families, namely symbolic and ordinary powers of ideals (see for example \cite{CutkoskyKurano} or \cite{CutkoskyHerzogTrung}). In contrast to Waldschmidt constant, which always exists for a graded family of monomial ideals, $\widehat{\reg}(\ibul)$ and $\widehat{\ri}(\ibul)$  may not exists.
The next two examples show that this situation can occur. These examples are of our interest also from other reasons. They show that dropping assumption about LBHR in Theorem \ref{th:ahf ahp polynomial monomial} may not affect on existence of $\ahp_{\ibul}$ and $\ahf_{\ibul}$.
\begin{example}
	\label{ex:ahp=1}
	Let $\Ibul\subset S(1)$. We define
	$$I_m =( x^2,xy^{2^{m}}) \text{ for } m\geq 1.$$
	Ideals $I_m$ fulfil condition 2) in Lemma \ref{th:borel_equivalent} directly from definition and therefore they are Borel-fixed.
	
	In order to show that $\ibul$ is a graded family observe that
	$$I_p I_q=(x^2,xy^{2^{p}})\cdot (x^2,xy^{2^{q}})=(x^4,x^3y^{2^{p}},x^3y^{2^{q}},x^2y^{2^p+2^{q}}) \subset (x^2,xy^{2^{q+q}})=I_{p+q}.$$ 
	By \cite[Proposition 2.11]{Green1} we see that $\reg(I_m) =\deg(xy^{2^{m}})= 2^{m}+1$.
	
	Fix a real number $t$. In this example we have
	\begin{equation*}
	\begin{split}
	M_{m,t}(\ibul)= \{a \in \R \; :\;a \geq 1 \text{ and } a\leq mt-2^{m} \}
	\cup \;\{ a \in \R\;:\;a \geq 2 \text{ and } a\leq mt \}.
	\end{split}
	\end{equation*}
	Directly from definitions of $\Gamma_{m,t}(\ibul)$ and $L_{m,t}(\ibul)$ we obtain
	$$\Gamma_{m,t}(\ibul)=\T_{mt}\setminus L_{m,t}(\ibul)=[0,mt] \setminus L_{m,t}(\ibul)= 
	\begin{cases}
	[0,1)\cup (mt-2^m,2) ,     & \text{ if  }mt-2^m\geq 1\\
	[0,2) , &  \text{ if }mt-2^m<1\\
	\end{cases}
	,$$
	thus
	\[ \widetilde{\HF}_{I_m}(t) =
	\begin{cases}
	1 ,     & \quad \text{ for  }t \geq \frac{2^m+1}{m}\\
	2 , & \quad \ \text{ for } t < \frac{2^m+1}{m}\\
	\end{cases},
	\]
	and as a consequence $\widetilde{\HF}_{I_m}(t)=\HP_{I_m}(t)=1$ for $t\geq \lfloor \frac{2^m+1}{m}\rfloor$. The limit$$\lim_{m \rightarrow \infty} \frac{\HP_{I_m}(mt)}{m}=\lim_{m \rightarrow \infty} \frac{1}{m}=0.$$
	From the $\widetilde{\HF}_{I_m}$ we can also read off that $\ri(I_m) =\lfloor \frac{2^m+1}{m}\rfloor$, which means that Hilbert regularity can not be linearly bounded and in addition
	$$\widehat{\reg}(\ibul) =\lim_{m \rightarrow \infty} \frac{2^{m}+1}{m}=+\infty=\lim_{m \rightarrow \infty} \frac{\lfloor \frac{2^m+1}{m}\rfloor}{m}= \widehat{\ri}(\ibul). $$
	
	All these considerations also show that in this example $\ahp_{\ibul}$ and $\ahf_{\ibul}$ exist although  LBHR does not hold.
\end{example}	
\begin{example}
	\label{ex:ahp_nonexistence}
	Take ideals $\jbul$ defined as in the Example \ref{ex:ahp=1} in polynomial ring $S(2)$. 
	
	Here the sets $M_{m,t}(\jbul)$ are
	\begin{equation*}
	\begin{split}
	M_{m,t}(\jbul)=  \{ (a,b) \in \R^2\; : \; 1\leq a \;\text{ and }\;  2^{m}\leq b\; \text{ and }\; a+b\leq mt-2^m \}\; \\
	\cup\;\{ (a,b) \in \R^2_{\geq 0}\; : \; 2\leq a\leq mt\; \text{ and }\; a+b \leq mt \}.
	\end{split}
	\end{equation*}
	In this case we have
	$$\Gamma_{m,t}(\jbul)=\T_{mt}\setminus L_{m,t}(\jbul)=\{(a,b) \in \R^2_{\geq 0}\; : \; a+b \leq mt\}\;\setminus L_{m,t}(\jbul)=$$
	$$
	=\begin{cases}
	\{(a,b) \in \R^2_{\geq 0}\; : 0\leq a<2 \text{ and } b<2^m-mt \} ,     & \text{ if  }mt\leq 2^m\\
	\{(a,b) \in \R^2_{\geq 0}\; : 0\leq a<2 \text{ and } b<2^m \} ,     & \text{ if  }mt> 2^m\\
	\end{cases}.$$
	Therefore
	\[ \widetilde{\HF}_{J_m}(t) =
	\begin{cases}
	2\lfloor t\rfloor+1 ,     & \quad \text{ for  }t \leq \frac{2^m}{m}\\
	\lfloor t\rfloor+1+2^m , & \quad \ \text{ for } t > \frac{2^m}{m}\\
	\end{cases}.
	\]
	
	We see from $\widetilde{\HF}_{J_m}$ that $\HP_{J_m}(t)=t+1+2^m$ and $ \ri(J_m) = \frac{2^m}{m}$, so LBHR does not hold.
	
	As a consequence, the limit $$\lim_{m \rightarrow \infty} \frac{\HP_{J_m}(mt)}{m^2}=\lim_{m \rightarrow \infty}\frac{\lfloor mt\rfloor+1+ 2^m}{m^2}$$ is infinity, so that the expression does not define any function. Again, as in Example \ref{ex:ahp=1} it can be easily seen that
	$$\widehat{\reg}(\jbul) = \widehat{\ri}(\jbul)=+\infty. $$
\end{example}
\begin{remark}
	Examples \ref{ex:ahp=1} and \ref{ex:ahp_nonexistence} can be generalized to more variables. It is enough to consider families of ideals $$ I_m=(x_0^2,x_0x_{1}^{2^m}), \subset S(n).$$
\end{remark}
The values of $\reg(I)$ and $\ri(I)$ can differ for Borel-fixed ideals. Take for instance the ideal $I=(x^5,x^4y,x^3y^3,x^2y^5,xy^7) \subset S(3)$. For this ideal $\ri(I)=6 < \reg(I)=\deg(xy^7)=8.$
It is an intriguing question however if their asymptotic counterparts may differ as well. What we can prove is the following.
\begin{theorem}
	Let $\Ibul$ be a graded sequence of Borel-fixed monomial ideals in $S(n)$. Assume that for each $m$, the ideals $I_m$ are Cohen–-Macaulay ideals and that $\widehat{\reg}(\ibul)$ exists. Then $\widehat{\reg}(\ibul)=\widehat{\ri}(\ibul)$.
\end{theorem}
\proof
It follows immediately from relation $\ri(I_m) = \reg(I_m)+\delta - n$, where $\delta$ is the projective dimension of $I_m$ (see \cite[Theorem 4.2]{Eis2}).
\endproof

\section{Limiting shape}
\label{sec:limiting shape}

Denote by $f_t\; : \; \R^n \rightarrow \R^{n+1}$ the function, which is defined as
$$f_t(x_0,\ldots,x_{n-1})=(x_0,\ldots,x_{n-1},t-\sum_{i=0}^{n-1} x_i).$$
Following the authors of articles \cite{Mayes, DSST, DST15}, we define the following geometrical object.

\begin{definition}
	\label{def:restricted and normal LimitingShape homogeneous}
	Let $\Ibul \subset S(n)$ be a graded system of homogeneous ideals. Fix a real number $t$. The restricted limiting shape of $\Ibul$ is
	$$\Delta(\ibul,t):= \overline{\bigcup_{m=1}^{\infty} \frac{L_{m,t}(\ibul)}{m}} \subseteq \R^n,$$
	while the limiting shape of $\Ibul$ is
	$$\Delta(\ibul):= \bigcup_{t \in \R} f_t\Big(\Delta(\ibul,t)\Big) \subseteq \R^{n+1},$$
	where the coordinates in $\R^{n}$ are $x_0,\ldots,x_{n-1}$ and in $\R^{n+1}$ they are $x_0,\ldots,x_{n-1},t$.
\end{definition}
\begin{definition}
	\label{def:restricted complementaryOfLimitingShape homogeneous}
	For the ideals $\Ibul \subset S(n)$ and a fixed number $t$ let
	$$\Gamma(\ibul,t):= \overline{\bigcap_{m=1}^{\infty}  \Big( \T_{mt}^n \setminus \frac{L_{m,t}(\ibul)}{m}\Big)}\subseteq \R^n,$$
	then the complement of the limiting shape of $\Ibul$ is
	$$\Gamma(\ibul):=\bigcup_{t \in \R} f_t\Big(\Gamma(\ibul,t)\Big)\subseteq \R^{n+1}.$$
\end{definition}

As the next step we want to show a connection between the volume of the complement of the limiting shape of $\Ibul$ and asymptotic Hilbert function. The same connection as authors of \cite[Theorem 11]{DST15} showed for homogeneous radical ideals. To do it we need the succeeding lemma

\begin{lemma}
	\label{lm:LsumVol homogeneous}
	Let $\Ibul \subset S(n)$ be a graded system of homogeneous ideals. For the sets $L_{m,t}(\ibul)$ there is
	$$ \vol\Big( \bigcup \frac{L_{m,t}(\ibul)}{m} \Big) = \vol\Big(  \overline{\bigcup \frac{L_{m,t}(\ibul)}{m}} \Big). $$
\end{lemma}
\proof
	Put $A:= \bigcup \frac{L_{m,t}(\ibul)}{m} $.
	
	We begin the proof by showing that $\partial(A) \subset \partial(\overline{A}) $, that is equivalently $\text{int}(\overline{A}) \subset \text{int}(A)$.
	
	Take $a=(a_0,\ldots,a_n) \in \text{int}(\overline{A}) \subset \overline{A} \subset \R^n$. By definition of sets $L_{m,t}$ we get that $a_0+\ldots+a_n \leq t$.
	From fact that $a\in \text{int}(\overline{A})$, we have that there exists a point $b=(b_0,\ldots,b_n) \in \overline{A}$ such that 
	$$b_0+\ldots+b_n < a_0+\ldots+a_n,$$
	and a ball $K(b,\epsilon)$, with $\epsilon >0$, such that $a \notin K(b,\epsilon)$. Since $b \in \overline{A}$ we may find
	$$x=(x_0,\ldots,x_n) \in K(b,\epsilon) \subset \text{int}(A).$$
	But then $x_0+\ldots+x_n < a_0+\ldots+a_n$ and that in fact means that $a \in \text{int}(A)$.
	
	From the definition of the sets $L_{m,t}(\ibul)$ we see that $A$ is bounded and from Theorem \ref{th:Lmt HomogeneousConvex} that $A$ is convex. Therefore $\vol(\partial(\overline{A}))=0$.
\endproof

\begin{theorem}
	\label{th:aHFvolGammaEqual homogeneous} 
	For a graded system of homogeneous ideals $\Ibul$, such that $\ini(I_m)$ is a Borel-fixed ideal for each $m$,
	$$\vol\Big(\Gamma(\ibul,t)\Big)= \ahf_{\ibul}(t). $$
\end{theorem}
\proof
	By Theorems \ref{th:ahf_homogeneous_ existence},\;\ref{th:gamma_vol homogeneous} and \ref{th:Lmt_sequence},we obtain
	$$\ahf_{\ibul}(t) = \lim_{m \rightarrow \infty}\vol\Big(\frac{\Gamma_{m,t}(\ibul)}{m}\Big) =\frac{t^n}{n!}- \lim_{m \rightarrow \infty}\vol\Big(\frac{L_{m,t}(\ibul)}{m}\Big)=\frac{t^n}{n!}- \vol\Big( \bigcup_m \frac{L_{m,t}(\ibul)}{m}\Big),$$
	applying Lemma \ref{lm:LsumVol homogeneous} we get
	$$\frac{t^n}{n!}- \vol\Big( \bigcup_m \frac{L_{m,t}(\ibul)}{m}\Big) = \frac{t^n}{n!}- \vol\Big(\overline{ \bigcup_m \frac{L_{m,t}(\ibul)}{m}}\Big) =$$
	$$=\vol\Big( \T_{mt}^n\setminus  \overline{\bigcup_{m=1}^{\infty} \frac{L_{m,t}(\ibul)}{m}}\Big)= \vol\Big( \T_{mt}^n\setminus \Delta(\ibul,t)\Big).$$
	
	Observe that from Lemma \ref{lm:LsumVol homogeneous} we have
	$$\vol\Big(\overline{\Delta(\ibul,t)}\setminus \text{int}(\Delta(\ibul,t))\Big)=0,$$
	so
	$$\vol\Big( \T_{mt}^n\setminus \Delta(\ibul,t)\Big)= \vol\Big(\overline{ \T_{mt}^n\setminus \Delta(\ibul,t)}\Big) =\vol \Big( \Gamma(\ibul,t)\Big).$$
\endproof

There is a nice geometrical connection between some already introduced asymptotic invariants and limiting shapes. This will be the subject of our next deliberations.
\begin{theorem}
	\label{th:WladschmidtAxis}
	Suppose that $\Ibul$ is a graded monomial Borel-fixed family of ideals. Then the value of $\widehat{\alpha}(\ibul)$ is the first coordinate of the point where for $t \gg 0$ sets $\Delta(\ibul,t)$ and $\Gamma(\ibul,t)$ meet the $x_0$-axis., i.e.
	$$\widehat{\alpha}(\ibul)=\min\{\alpha_0\;:\;(\alpha_0,0,\ldots,0)\in\Gamma(\ibul,t)\cap \{x_1=\ldots=x_n=0\}\subseteq \R^{n+1}\}. $$
\end{theorem}
\proof
	 Assume that $\{g_{m,1},\ldots,g_{m,s_m}\}$ are the minimal sets of generators of $I_m$ with $\deg(g_{m,1})\geq\ldots\geq \deg(g_{m,s_m})$. Since every ideal $I_m$ is a Borel-fixed ideal, applying Lemma \ref{th:borel_equivalent} we obtain that
	$$x_0^{\deg(g_{m,1})},\;x_0^{\deg(g_{m,2})}, \ldots, x_0^{\deg(g_{m,s_m})} \in I_m.$$
	The Waldschmidt constant depends on the values of $\alpha(I_m)$, which is equal to
	$\alpha(I_m)=\deg(g_{m,s_m})$.
	This implies
	$$ (\deg(g_{m,1}),0,\ldots,0),\; \ldots, (\deg(g_{m,s_m}),0,\ldots,0)\in L_{m,t}(\ibul)$$
	and $$ (\deg(g_{m,s_m}),0,\ldots,0)\in \overline{\Gamma_{m,t}(\ibul) },$$
	for $m$ and $t$ big enough.
	
	Since
	$$\widehat{\alpha}(\ibul)= \inf_m \frac{\deg(g_{s_m})}{m},$$
	we obtain
	$$\Big(\widehat{\alpha}(\ibul),0,\ldots,0\Big) \in \bigcup\Big( \overline{\frac{L_{m,t}(\ibul)}{m}}\Big) \subseteq \Delta(\ibul,t),$$ and
	$$\Big(\widehat{\alpha}(\ibul),0,\ldots,0\Big) \in \bigcap\Big( \overline{\frac{\Gamma_{m,t}(\ibul)}{m}}\Big)\subseteq \Gamma(\ibul,t).$$
\endproof
Before we set similar description for asymptotic regularity we need to present a specific example to show that LBHR is not enough to guarantee the existence of this invariant. In some cases the bound for the regularity has to be more tight.
\begin{example}
	\label{ex:reg i 2 pkt skupienia}
	Let $a,b,d$ be fixed integers with $a<b$ and $d\geq2$. Consider the graded family $\Ibul$ of ideals defined as
	$$I_{dk-d+1}=I_{d(k-1)+1}=(x^{ak}), $$
	$$I_{dk-d+2}=I_{d(k-1)+2}=(x^{ak+1},x^{ak}y^{bk}), $$
	$$\vdots$$
	$$I_{dk}=(x^{ak+1},x^{ak}y^{bk}),$$
	for all $k \in \{1,2,\ldots \}$.
	We may check by hand that all these ideals are Borel-fixed. We want to show that the condition $I_p \cdot I_q \subset I_{p+q}$ holds for all $p,q\in \{1,2,\ldots\}$. Thus we need to consider three cases:
	\begin{itemize}
		\item [1)] $I_{d(p-1)+1}\cdot I_{d(q-1)+1} \subset I_{d((p+q)-2)+2}$ if, and only if
		$x^{ap}\cdot x^{aq}=x^{a(p+q)}$ belongs to the right side of the condition. By definition $x^{a(p+q-1)+1} \in I_{d((p+q)-2)+2}$, which gives the inclusion.
		\item [2)] In case $I_{d(p-1)+1}\cdot I_{d(q-1)+r} \subset I_{d((p+q)-2)+r+1}$ with $r \in \{2,3\ldots\,d\}$, we check if
		$x^{ap}\cdot x^{aq+1}=x^{a(p+q)+1}$ and $x^{ap}\cdot x^{aq}y^{bq}=x^{a(p+q)}y^{bq}$
		belongs to $I_{d((p+q)-2)+r+1}$. In this case we use monomials
		$$x^{a(p+q-1)+1}\quad \text{ if }\quad r\leq d-1, $$
		$$x^{a(p+q-1)}\quad \text{ if }\quad r=d-1 ,$$
		coming from $I_{d((p+q)-2)+r+1}$ to show divisibility in both subcases.
		\item [3)] In the case $I_{d(p-1)+r_1}\cdot I_{d(q-1)+r_2} \subset I_{d((p+q)-2)+r_1+r_2}$ with $r_1,r_2  \in \{2,3\ldots\,d\}$, we have
		$x^{ap+1}\cdot x^{aq+1}=x^{a(p+q)+2}$ and $x^{ap+1}\cdot x^{aq}y^{bq}=x^{a(p+q)+1}y^{bq}$. These monomials are divisible by
		\[
		\begin{cases}
		x^{a(p+q-1)+1} ,     & \quad \text{ for  }r_1+r_2 \leq d,\\
		x^{a(p+q)} , & \quad  \text{ for }r_1+r_2 = d+1,\\
		x^{a(p+q)+1}, & \quad \text{ for } d+1<r_1+r_2 \leq 2d.\\
		\end{cases}
		.
		\]
		For monomial $x^{ap}y^{bp}\cdot x^{aq}y^{bq}=x^{a(p+q)}y^{b(p+q)}$ we use
		\[
		\begin{cases}
		x^{a(p+q-1)} y^{b(p+q-1)},     & \quad \text{ when  }r_1+r_2 \leq d,\\
		x^{a(p+q)} , & \quad  \text{ for }r_1+r_2 = d+1,\\
		x^{a(p+q)} y^{b(p+q)}, & \quad \text{ when } d+1<r_1+r_2 \leq 2d.\\
		\end{cases}
		.
		\]to show divisibility. All these monomials are in $I_{d((p+q)-2)+r_1+r_2}$ by definition of $\ibul$.
	\end{itemize}
	Simply calculations show that asymptotic regularity does not exist, indeed
	$$ \lim_{k \rightarrow \infty}\frac{\reg(I_{d(k-1)+1})}{d(k-1)+1}= \lim_{k \rightarrow \infty}\frac{ak}{d(k-1)+1}=\frac{a}{d}, $$
	$$ \lim_{k \rightarrow \infty}\frac{\reg(I_{dk})}{dk}= \lim_{k \rightarrow \infty}\frac{(a+b)k}{dk}=\frac{a+b}{d} \neq \frac{a}{d}.$$
\end{example}
These considerations lead to the following.
\begin{lemma}
	\label{lm:aregExisting homogeneous} 
	Suppose that $\Ibul$ is a graded homogeneous family of ideals such that $\ini(I_m)$ is a Borel-fixed ideal for all $m$ and
	$$am+b \leq \reg(I_m) \leq am+c, $$
	for some $a \geq 0$, $0\leq b \leq c$ and $m$ big enough. Then $\widehat{\reg}(\ibul)$ always exists and is equal $a$.
\end{lemma}
\proof
	The proof of this lemma is obvious. It is worth pointing out here that the previous example shows that this lemma may not hold if we drop any of the restrictions in the statement.
\endproof
Now we define a new quantity for the sake of the next theorem.
\begin{definition}
	For a graded family of homogeneous ideals $\ibul$ let
	$$a(\ibul):=\sup \{|x|, \text{ such that } x=(x_0,x_1,\ldots,x_n) \text{ is an extremal point in }\Delta(\ibul)\} \in \R\cup \{\infty\}.$$
\end{definition}
\begin{theorem}
	\label{th:aregPointGin homogeneous} 
	Let $\ibul$ be a graded family of homogeneous ideals such that $\ini(I_m)$ is Borel-fixed ideal for all $m$. If $\widehat{\reg}(\ibul)$ exists, then
	$$\widehat{\reg}(\ibul)=a(\ibul).$$ 
\end{theorem}
\proof
	Assume that $\widehat{\reg}(\ibul)$ exists. Fix $t \in \R$ and observe that every extremal point in $\conv(L_{m,t}(\ibul))$ comes from one of the generators of the minimal set of generators of $I_m$, or this point can be obtained by condition 2) from Lemma \ref{th:borel_equivalent}. Extremal point of $\Delta(\ibul,t)$ is a limit of extremal points of  $\conv(L_{m,t}(\ibul))$ (see \cite[Theorem 2]{Jerison}). Therefore for every $t$, from convexity of $\overline{\bigcup\frac{L_{m,t}(\ibul)}{m}}$, every extremal point in $\Delta(\ibul,t)$ comes from one of the generators of the minimal set of generators of $I_m$.
	Borel-fixed ideals are monomial ideals, thus from \cite[Proposition 2.11]{Green1} we know that for any $m$, $\reg(I_m)$ is equal to the maximal degree of the generator in the minimal set of generators. We conclude that $\widehat{\reg}(\ibul)$ must be equal to $a(\ibul)$.
\endproof
Next examples show that for a graded sequence of ideals it is possible to construct a family of ideals for which the Waldschmidt constant and asymptotic regularity are irrational numbers. We recall here that we can find in the literature an example of irrational asymptotic regularity calculated for ordinary powers of the ideal sheaf on $\P^n$. See \cite{Cutkosky} for more details.
\begin{example}
	\label{ex:Waldschmidt constant irrational}
	Let $\ibul=\{I_m\}_m \subset \mathbb{K}[x]$ be a family of ideals defined as $I_m=(x^{\lceil m \pi\rceil})$.
	
	First, we show that it is indeed a graded sequence. To see that, take $x^{\lceil p \pi\rceil} \in I_p$ and $x^{\lceil q \pi\rceil} \in I_q$, then
	$$ x^{\lceil p \pi\rceil} \cdot x^{\lceil q \pi\rceil} =  x^{\lceil q \pi\rceil + \lceil p \pi\rceil} \in I_{p+q},$$
	since we know that inequality $\lceil a+b\rceil \leq \lceil a \rceil + \lceil b \rceil$ holds for any numbers $a,b$.
	
	Finally, we calculate that
	$$ \widehat{\alpha}(\ibul) =\lim_{m \rightarrow \infty} \frac{\alpha(I_m)}{m}=\lim_{m \rightarrow \infty}\frac{\lceil m \pi\rceil}{m}=\pi.$$
\end{example}
\begin{example}
	\label{ex:Waldschmidt areg irrational}
	Pick two positive numbers $q_1 \leq q_2$. Define the ideals $I_m=(x^{(a,b)}) \subset S(1)$, where $a,b$ are all positive integers such that the inequality
	$aq_2 + bq_1 \geq mq_1 q_2$ is fulfilled.
	
	As in Example \ref{ex:Waldschmidt constant irrational}, we start by showing that $\ibul=\{I_m\}_m$ is a graded sequence.
	Let $x^{(a_k,b_k)} \in I_k$ and $x^{(a_l,b_l)} \in I_l$, then by definition
	$$a_kq_2 + b_kq_1 \geq kq_1 q_2 \;\text{ and }\;\ a_lq_2 + b_lq_1 \geq lq_1 q_2,$$
	from which we get
	\begin{equation*}
	\label{eq:Waldschmidt areg irrational}
	(a_k+a_l)q_2 + (b_k+b_l)q_1 \geq (k+l)q_1 q_2.
	\end{equation*}
	This implies
	$$x^{(a_k,b_k)} \cdot x^{(a_l,b_l)} = x^{(a_k+a_l,b_k+b_l)} \in I_{k+l}.$$
	
	Now we show that $\ibul$ is a Borel-fixed family of ideals. Pick any $x^{(a,b)} \in I_m$. We want to check if $x^{(a+1,b-1)} \in I_m$. If so, then from Lemma \ref{th:borel_equivalent} we obtain that $I_m$ is a Borel-fixed ideal.
	
	We have
	$$ (a+1)q_2 + (b-1)q_1 \geq aq_2 + bq_1 +q_2 - q_1 \geq mq_1 q_2 + q_2 - q_1 \geq mq_1 q_2,$$
	which leads to desired membership.
	
	From Theorems \ref{th:WladschmidtAxis} and \ref{th:aregPointGin homogeneous}  we obtain
	$$ \widehat{\alpha}(\ibul) =q_1\; \text{ and }\;\widehat{\reg}(\ibul)=q_2.$$
\end{example}
\begin{remark}
	The number $\pi$ in Example \ref{ex:Waldschmidt constant irrational} can be easily replaced by any other irrational number. The result remains unchanged. We also see that in Example \ref{ex:Waldschmidt areg irrational} the numbers $q_1$ and $q_2$ can be irrational. In both examples we can increase the number of variables, ,i.e., considering these ideals in any ring $S(n)$, $n\geq 2$. 
\end{remark}
It is possible to find a graded sequence of ideals for which we can construct a limiting shape with finite many hyperplane segments forming its boundary. All we have to do is to proceed in similar way as in Example \ref{ex:Waldschmidt areg irrational}. In the next example we can see how it works for ideals in $S(1)$. The more general case is tedious and therefore here omitted.  
\begin{example}
	\label{ex:anyDescent}
	In this example we create a family of ideals in $S(1)$ such that the boundary of limiting shape is as indicated on Figure \ref{fig:anyDescent}. In order to make it, every ideal $I_m$ has to contains all monomials which are in the convex set defined by points $(0,mt_n),\ldots,(ms_0,0)$ (see Figure \ref{fig:anyDescentM}).
	
	\begin{minipage}{0.49\textwidth}
		\begin{figure}[H]
			\begin{tikzpicture}[line cap=round,line join=round,>=triangle 45,x=1.0cm,y=1.0cm,scale=0.7]
			\clip(-0.34,-0.5) rectangle (7.5,8.24);
			\draw [->,line width=1.pt] (0.,0.) -- (0.,8.);
			\draw [->,line width=1.pt] (0.,0.) -- (7.,0.);
			\draw [line width=1.pt] (0.,6.44)-- (0.56,3.38);
			\draw [line width=1.pt] (0.56,3.38)-- (1.12,2.04);
			\draw [line width=1.pt,dotted] (1.12,2.04)-- (2.32,0.86);
			\draw [line width=1.pt] (2.32,0.86)-- (4.8,0.);
			\begin{scriptsize}
			\draw (0.12,7.) node[anchor=north west] {$(0,t_n)$};
			\draw (0.66,3.96) node[anchor=north west] {$(s_{n-1},t_{n-1})$};
			\draw (1.26,2.7) node[anchor=north west] {$(s_{n-2},t_{n-2})$};
			\draw (2.36,1.64) node[anchor=north west] {$(s_1,t_1)$};
			\draw (4.62,0.88) node[anchor=north west] {$(s_0,0)$};
			\draw [fill=black] (0.,6.44) circle (2.5pt);
			\draw [fill=black] (0.56,3.38) circle (2.5pt);
			\draw [fill=black] (1.12,2.04) circle (2.5pt);
			\draw [fill=black] (2.32,0.86) circle (2.5pt);
			\draw [fill=black] (4.8,0.) circle (2.5pt);
			\end{scriptsize}
			\end{tikzpicture}
			\caption{.  The boundary of the limiting shape of $\Delta(\ibul)$.}
			\label{fig:anyDescent}
		\end{figure}
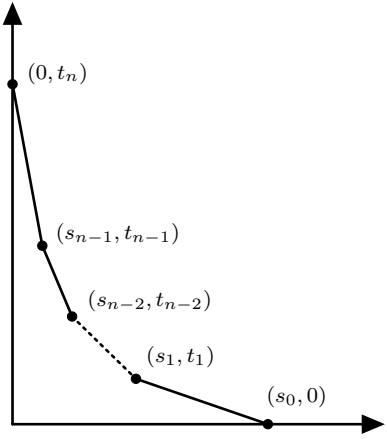
	\end{minipage}
	\begin{minipage}{0.49\textwidth}
			\begin{figure}[H]
			\begin{tikzpicture}[line cap=round,line join=round,>=triangle 45,x=1.0cm,y=1.0cm,scale=0.7]
			\clip(-0.34,-0.5) rectangle (7.5,8.24);
			\draw [->,line width=1.pt] (0.,0.) -- (0.,8.);
			\draw [->,line width=1.pt] (0.,0.) -- (7.,0.);
			\draw [line width=1.pt] (0.,6.44)-- (0.56,3.38);
			\draw [line width=1.pt] (0.56,3.38)-- (1.12,2.04);
			\draw [line width=1.pt,dotted] (1.12,2.04)-- (2.32,0.86);
			\draw [line width=1.pt] (2.32,0.86)-- (4.8,0.);
			\begin{scriptsize}
			\draw (0.12,7.) node[anchor=north west] {$(0,mt_n)$};
			\draw (0.66,3.96) node[anchor=north west] {$(ms_{n-1},mt_{n-1})$};
			\draw (1.26,2.7) node[anchor=north west] {$(ms_{n-2},mt_{n-2})$};
			\draw (2.36,1.64) node[anchor=north west] {$(ms_1,mt_1)$};
			\draw (4.62,0.88) node[anchor=north west] {$(ms_0,0)$};
			\draw [fill=black] (0.,6.44) circle (2.5pt);
			\draw [fill=black] (0.56,3.38) circle (2.5pt);
			\draw [fill=black] (1.12,2.04) circle (2.5pt);
			\draw [fill=black] (2.32,0.86) circle (2.5pt);
			\draw [fill=black] (4.8,0.) circle (2.5pt);
			\end{scriptsize}
			\end{tikzpicture}
			\caption{.  The boundary of the ideals $I_m$.}
			\label{fig:anyDescentM}
		\end{figure}
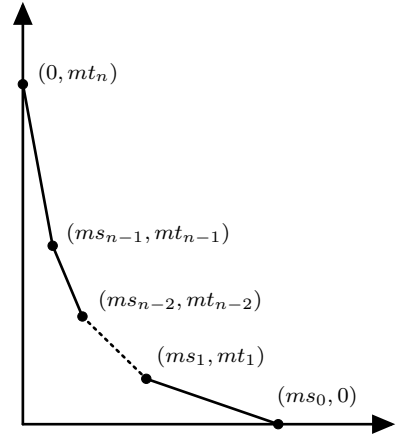
	\end{minipage}
It can be easily obtained by considering the equations of lines which contains each segment. Thus
\begin{equation}
	\begin{vmatrix}
		ms_{i+1} & mt_{i+1} & 1  \\ 
		ms_{i} & ms_{i} & 1  \\ 
		x & y & 1   \notag
	\end{vmatrix} \geq 0
	\iff x(t_{i+1}-t_i)+y(s_i-s_{i+1})\geq m(s_it_{i+1}-s_{i+1}t_i), 
\end{equation}
for $i\in \{0,\ldots,n-1\} $, are desired equations. 

The condition 2) from Lemma \ref{th:borel_equivalent} and the convexity of $\Delta(\ibul)$ are fulfilled if
\begin{equation}
\frac{t_{n}-t_{n-1}}{s_{n}-s_{n-1}}< \ldots < \frac{t_{i+1}-t_i}{s_{i+1}-s_i}< \ldots < \frac{t_{1}-t_0}{s_{1}-s_0} \leq -1.
\label{eq:numbers}
\end{equation}
We are ready to set the definition of ideals $I_m$. 

Pick positive numbers $t_i$ and $s_i$ such that condition (\ref{eq:numbers}) is fulfilled. Let  $I_m=(x^{(a,b)}) \subset S(1)$, where $a,b$ are positive integers satisfying 
$a(t_{i+1}-t_i)+b(s_i-s_{i+1})\geq m(s_it_{i+1}-s_{i+1}t_i)$ for all $i \in \{0,\ldots,n-1\} $. Then $\Delta(\ibul)$ is the convex set bounded by points $(s_0,0),(s_{n-1},t_{n-1}),\ldots,(0,t_n)$ and $x,y$-axis.

To check that this definition gives rise to graded sequence of ideal take $x^{(p_1,p_2)} \in I_p$ and $x^{(q_1,q_2)} \in I_q$ sum up the inequalities
$p_1(t_{i+1}-t_i)+p_2(s_i-s_{i+1})\geq p(s_it_{i+1}-s_{i+1}t_i) $ and 
$q_1(t_{i+1}-t_i)+q_2(s_i-s_{i+1})\geq q(s_it_{i+1}-s_{i+1}t_i) $
for every number $i$.

Finally, let $x^{(a,b)} \in I_m$ and observe that for every $i$ we have
$(a+1)(t_{i+1}-t_i)+(b-1)(s_i-s_{i+1})\geq m(s_it_{i+1}-s_{i+1}t_i) + (t_{i+1}-t_i)-(s_i-s_{i+1})\geq m(s_it_{i+1}-s_{i+1}t_i) $ by (\ref{eq:numbers}), from which we deduce that $x^{(a+1,b-1)} \in I_m$.
\end{example}
\begin{remark}
	It is still an open problem to find a graded sequence of ideals formed by symbolic or ordinary powers for which the breaking points of limiting shape have irrational coordinates. But it is possible to construct examples with any, but finitely many such points. See Corollary \ref{cor:n break points gin}.
\end{remark}

\section{Asymptotic first difference Hilbert function}
\label{sec:connection aHF daHF}
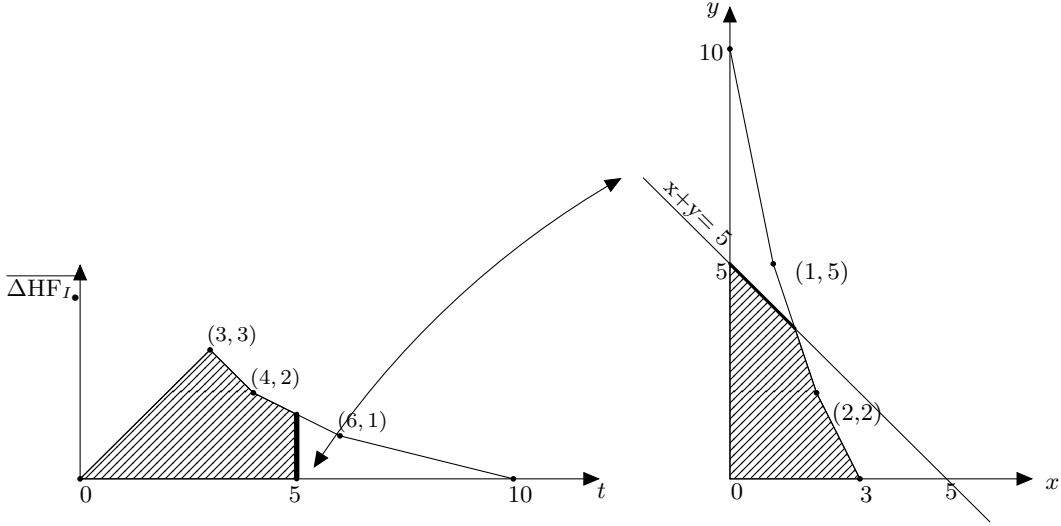
\begin{figure}
	\centering
	\begin{tikzpicture}[line cap=round,line join=round,>=triangle 45,x=1.0cm,y=1.0cm,scale=0.57]
	\clip(-1.9024669415683841,-1.5835970902371785) rectangle (23.041629808095454,11.29278657043657);
	\fill[fill=black,pattern=north east lines,pattern color=black] (18.,0.) -- (17.,2.) -- (16.5,3.5) -- (15.,5.) -- (15.,6.277519200586862E-5) -- cycle;
	\fill[fill=black,pattern=north east lines,pattern color=black] (0.,0.) -- (3.,3.) -- (4.,2.) -- (5.,1.5) -- (5.,0.) -- cycle;
	\draw [->] (0.,0.) -- (12.,0.);
	\draw (0.,0.)-- (3.,3.);
	\draw (3.,3.)-- (4.,2.);
	\draw (6.,1.)-- (4.,2.);
	\draw (6.,1.)-- (10.,0.);
	\draw (15.,10.)-- (16.,5.);
	\draw (16.,5.)-- (17.,2.);
	\draw (17.,2.)-- (18.,0.);
	\draw [line width=2.pt] (5.,0.)-- (5.,1.5);
	\draw [->] (15.,0.) -- (15.,11.);
	\draw (18.,0.)-- (17.,2.);
	\draw (17.,2.)-- (16.5,3.5);
	\draw [line width=1.2pt] (16.5,3.5)-- (15.,5.);
	\draw (15.,5.)-- (15.,6.277519200586862E-5);
	\draw (15.,6.277519200586862E-5)-- (18.,0.);
	\draw (0.,0.)-- (3.,3.);
	\draw (3.,3.)-- (4.,2.);
	\draw (4.,2.)-- (5.,1.5);
	\draw (5.,1.5)-- (5.,0.);
	\draw (5.,0.)-- (0.,0.);
	\draw (15.,5.)-- (16.5,3.5);
	\draw (13.,7.)-- (21.,-1.);
	\begin{footnotesize}
	\draw (12.9,7.5) node[anchor=north west] {$\rotatebox{-45.0}{ \text{ x+y= 5 }  }$};
	\draw (-0.2106419566421359,0.028872007015976982) node[anchor=north west] {$0$};
	\draw (4.613718273929303,0.028872007015976982) node[anchor=north west] {$5$};
	\draw (9.674494683040276,0.02897935907046201) node[anchor=north west] {$10$};
	\draw (17.80694752136056,0.02897935907046201) node[anchor=north west] {$3$};
	\draw (14.815700500369191,0.1) node[anchor=north west] {$0$};
	\draw (11.730977009971841,0.09897935907046201) node[anchor=north west] {$t$};
	\end{footnotesize}
	\draw [->] (15.,0.) -- (22.,0.);
	\draw [->] (0.,0.) -- (0.,5.);
	\begin{footnotesize}
	\draw (-1.9,4.959755768181424) node[anchor=north west] {$\overline{\dHF_{\ibul}}$};
	\draw (22.083495996684157,0.21582494582793704) node[anchor=north west] {$x$};
	\draw (14.25484168393331,11.269417453085076) node[anchor=north west] {$y$};
	\draw (16.28795489351338,5.263554293750859) node[anchor=north west] {$(1,5)$};
	\draw (17.1526122355187,1.9918778645415578) node[anchor=north west] {(2,2)};
	\draw (14.0810344923032,10.3) node[anchor=north west] {$10$};
	\draw (14.465163740096765,5.2) node[anchor=north west] {$5$};
	\draw (19.769953378886147,0.1) node[anchor=north west] {$5$};
	\end{footnotesize}
	\draw [<->] (12.5,7) arc (120:147:21);
	\begin{scriptsize}
	\draw [fill=black] (0.,0.) circle (1.5pt);
	\draw [fill=black] (3.,3.) circle (1.5pt);
	\draw[color=black] (3.5050477022455793,3.323917553576773) node {$(3, 3)$};
	\draw [fill=black] (4.,2.) circle (1.5pt);
	\draw[color=black] (4.5099197483598665,2.3190455074624876) node {$(4, 2)$};
	\draw [fill=black] (6.,1.) circle (1.5pt);
	\draw[color=black] (6.519663840588443,1.3375425786996975) node {$(6, 1)$};
	\draw [fill=black] (10.,0.) circle (1.5pt);
	\draw [fill=black] (18.,0.) circle (1.5pt);
	\draw [fill=black] (17.,2.) circle (1.5pt);
	\draw [fill=black] (16.,5.) circle (1.5pt);
	\draw [fill=black] (15.,10.) circle (1.5pt);
	\draw [fill=black] (5.,0.) circle (0.5pt);
	\draw [fill=black] (5.,1.5) circle (0.5pt);
	\draw [fill=black] (16.5,3.5) circle (0.5pt);
	\draw [fill=black] (15.,5.) circle (0.5pt);
	\end{scriptsize}
	\end{tikzpicture}
	\caption{. A connection between $\overline{\dHF_{\ibul}}(t)$ and $\vol(\Gamma(\ibul,t))$ for $t=5$.}\label{fig:dhf & gin}
\end{figure}

For any function $f: \mathbb{Z} \longrightarrow \mathbb{Z}$ we write $\Delta f$ for the first difference function $ \Delta f(t): = f(t)-f(t-1).$
In this section we want to prove a nice connection between the theory for Hilbert functions and limiting shapes of graded systems of homogeneous ideals. The main theorem in this part is the following.
\begin{theorem}
	\label{th:adHF aHF connection}
	Let $\Ibul$ be a graded sequence of homogeneous ideals for which $\ahf_{\ibul}$ exists. Assume that for all $t\in\R$ we have that the limit
	$$\overline{\Delta\HF_{\ibul} }(t):=\lim_{m \rightarrow \infty} \frac{\Delta \widetilde{\HF}_{I_m}(mt)}{m^{n-1}}$$
	exists and is non-negative.
	For any real positive number $t_0$, denote by $V(t_0)$ the volume of the body between
	the graphs of $\overline{\dHF_{\ibul}}$, $x=t_0$ and the $t$--axis. Then
	$$V(t_0)=\ahf_{\ibul}(t_0).$$
\end{theorem}
\proof
For the given function $\overline{\dHF_{\ibul}}$ and $t_0 \in \R_{\geq 0}$ we divide the interval $[0,\lfloor t_0 \rfloor]$
into $s$ subintervals $[i \frac{1}{s},(i+1)\frac{1}{s}]$, for $i=0,\ldots,\lfloor st_0 \rfloor -1$, and $[\lfloor st_0 \rfloor\frac{1}{s}, t_0]$. Then
$$V(t_0)=\lim_{s \rightarrow \infty} \Bigg(\sum_{i=0}^{\lfloor st_0\rfloor} \overline{\Delta\HF_{\ibul}}\Big(\frac{i}{s}\Big)\frac{1}{s} + \Big(t_0-\frac{1}{s} \lfloor s t_0 \rfloor \Big)\overline{\Delta\HF_{\ibul}}(t_0)\Bigg)=$$
$$\lim_{s \rightarrow \infty} \lim_{m \rightarrow \infty} \frac{ \sum_{i=0}^{\lfloor st_0\rfloor}\Delta\widetilde{\HF}_{I_m} (\frac{i}{s}m)\frac{1}{s}}{m^{n-1} }$$
because $\Big(t_0-\frac{1}{s} \lfloor s t_0 \rfloor \Big)\overline{\Delta\HF_{\ibul}}(t_0) \rightarrow 0$ in the same time, when $ s \rightarrow \infty$.
Thus
$$V(t_0)=\lim_{s \rightarrow \infty} \lim_{m \rightarrow \infty} \frac{\Delta\widetilde{\HF}_{I_m}(m\frac{\lfloor st_0\rfloor}{s})\frac{1}{s} - \Delta\widetilde{\HF}_{I_m}(0)\frac{1}{s}}{m^{n-1}}$$
Since we assume that this limit exists, we use the so-called ``diagonal sequence trick'' (see e.g. \cite[Theorem 1.24]{ReedSimon}) and put $s=m$. Then
$$ V(t_0)=\lim_{m \rightarrow \infty}\frac{\widetilde{\HF}_{I_m}(\lfloor mt_0 \rfloor)}{m^n} =\lim_{m \rightarrow
	\infty}\frac{\widetilde{\HF}_{I_m}(mt_0)}{m^n}=\ahf_{\ibul}(t_0).$$
\endproof
Using the previous theorem and Theorem \ref{th:aHFvolGammaEqual homogeneous} we may reproduce the complement of the limiting shape in any cases for which we know $\overline{\Delta\HF_{\ibul}}$. The next corollary shows how these two theorems work in $\P^2$. An example of its application is indicated in the Figure \ref{fig:dhf & gin} (which presents just an idea, it is not associated to any specific graded system of ideals). The volumes of the hatched areas are equal.
By comparing the equal volumes we may obtain the following description for $\P^2$.
\begin{corollary}
	\label{cor:reading gin from dHF graph in P^2}
	Consider a graded family of homogeneous ideals $\ibul$ such that $I_m \subset S(2)$ for all $m$, and for which the assumptions of Theorem \ref{th:adHF aHF connection} are fulfilled. Then, every point $(x,y)$, with $x\leq t$, taken from the graph of $\overline{\dHF_{\ibul}}$ gives the point $(y,x-y)$, which belongs to the boundary of $\Gamma(\ibul,t)$.
\end{corollary}
\proof

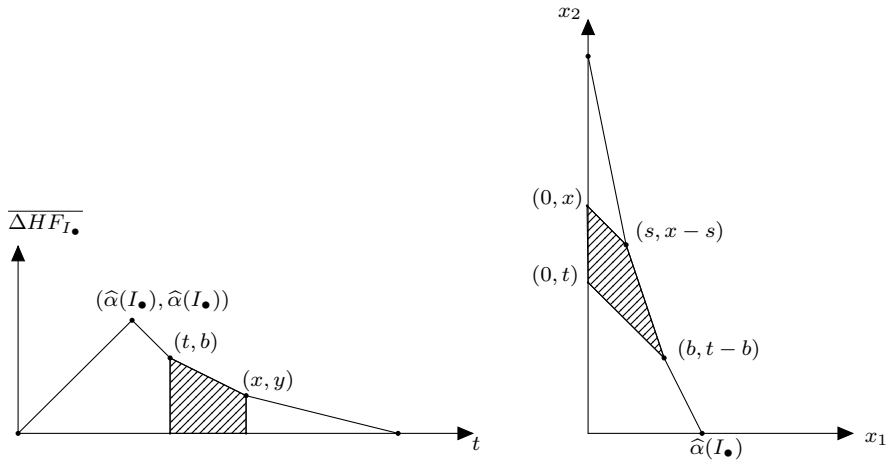
\begin{figure}[H]
	\begin{tikzpicture}[line cap=round,line join=round,>=triangle 45,x=1.0cm,y=1.0cm,scale=0.5]
	\clip(-0.5798740107957499,-0.9269248926601652) rectangle (22.838318487137446,11.283438923495975);
	\fill[fill=black,pattern=north east lines,pattern color=black] (4.,0.) -- (6.,0.) -- (6.,1.) -- (4.,2.) -- cycle;
	\fill[fill=black,pattern=north east lines,pattern color=black] (15.,4.) -- (17.,2.) -- (16.,5.) -- (14.97227358662421,6.039408989820496) -- cycle;
	\draw [->] (0.,0.) -- (12.,0.);
	\draw (0.,0.)-- (3.,3.);
	\draw (3.,3.)-- (4.,2.);
	\draw (6.,1.)-- (4.,2.);
	\draw (6.,1.)-- (10.,0.);
	\draw (15.,10.)-- (16.,5.);
	\draw (16.,5.)-- (17.,2.);
	\draw (17.,2.)-- (18.,0.);
	\draw [->] (15.,0.) -- (15.,11.);
	\draw [->] (15.,0.) -- (22.,0.);
	\draw [->] (0.,0.) -- (0.,5.);
	\draw (4.,2.)-- (4.,0.);
	\draw (6.,1.)-- (6.,0.);
	\draw (4.,0.)-- (6.,0.);
	\draw (6.,0.)-- (6.,1.);
	\draw (6.,1.)-- (4.,2.);
	\draw (4.,2.)-- (4.,0.);
	\draw (14.97227358662421,6.039408989820496)-- (16.,5.);
	\draw (15.,4.)-- (17.,2.);
	\draw (15.,4.)-- (17.,2.);
	\draw (17.,2.)-- (16.,5.);
	\draw (16.,5.)-- (14.97227358662421,6.039408989820496);
	\draw (14.97227358662421,6.039408989820496)-- (15.,4.);
	\begin{scriptsize}
	\draw (-0.4770498944491716,6.142233106167074) node[anchor=north west] {$\overline{\Delta HF_{I_{\bullet}}}$};
	\draw (22.067137614538108,0.22984641623883761) node[anchor=north west] {$x_1$};
	\draw (14,11.509144952582619) node[anchor=north west] {$x_2$};
	\draw (1.8107866942621944,4) node[anchor=north west] {$(\widehat{\alpha}(I_{\bullet}),\widehat{\alpha}(I_{\bullet}))$};
	\draw (3.8672690211937595,2.9) node[anchor=north west] {$(t,b)$};
	\draw (5.672,1.9) node[anchor=north west] {$(x,y)$};
	\draw (17.465758408028734,0.07561024171897057) node[anchor=north west] {$\widehat{\alpha}(I_{\bullet})$};
	\draw (11.733313921706996,0.12702229989225958) node[anchor=north west] {$t$};
	\draw (17.131580029902352,2.7233312376433547) node[anchor=north west] {$(b,t-b)$};
	\draw (13.3,4.676989448228337) node[anchor=north west] {$(0,t)$};
	\draw (16.026220779176636,5.83376075712734) node[anchor=north west] {$(s,x-s)$};
	\draw (13.3,6.707765746073253) node[anchor=north west] {$(0,x)$};
	\draw [fill=black] (0.,0.) circle (1.5pt);
	\draw [fill=black] (3.,3.) circle (1.5pt);
	\draw [fill=black] (4.,2.) circle (1.5pt);
	\draw [fill=black] (6.,1.) circle (1.5pt);
	\draw [fill=black] (10.,0.) circle (1.5pt);
	\draw [fill=black] (18.,0.) circle (1.5pt);
	\draw [fill=black] (17.,2.) circle (1.5pt);
	\draw [fill=black] (16.,5.) circle (1.5pt);
	\draw [fill=black] (15.,10.) circle (1.5pt);
	\draw [fill=black] (4.,0.) circle (0.5pt);
	\draw [fill=black] (6.,0.) circle (0.5pt);
	\draw [fill=black] (15.,4.) circle (0.5pt);
	\draw [fill=black] (14.97227358662421,6.039408989820496) circle (0.5pt);
	\end{scriptsize}
	\end{tikzpicture}
	\caption{.  The volume of hatched area presents the numbers $V(x)-V(t)$ and $P(x)-P(t)$.}\label{fig:dhf & gin area}
\end{figure}

Denote by $V(t_0)$ the volume of the body between the graphs of $\overline{\dHF_{\ibul}}$, $x=t_0$ and $t$--axis. By $P(t_0)$ we denote the volume of the area cut from $\Gamma(\ibul,t)$ by the "simplex" boundary by the coordinate axes and the line $x+y=t_0$.
Observe that for every point from the graph $\overline{\dHF_{\ibul}}$ of the form $(t,t)$, and with assumption $t\leq \widehat{\alpha}(\ibul)$, we have
$V(t)=\frac{1}{2} t^2=P(t). $
Now assume that for some $t \geq \widehat{\alpha}(\ibul)$ we have $V(t)=P(t)$. Let $b=\overline{\dHF_{\ibul}}(t)$. Take any point $x$, such that $t<x$. For such point, put $y=\overline{\dHF_{\ibul}}(x)$. Figure \ref{fig:dhf & gin area} presents the situation already described, with the point $(s,x-s)$, which coordinates we want to determine.
From Theorems \ref{th:adHF aHF connection} and \ref{th:aHFvolGammaEqual homogeneous} we get
$V(x)-V(t)=P(x)-P(t),$
which means that
$$V(x)-V(t)=\frac{y+b}{2}(x-t)=\frac{1}{2} \Bigg(
\begin{vmatrix}
-s  & s \\
0 & x-t  \notag
\end{vmatrix}
+
\begin{vmatrix}
s-b & x-s-t+b \\
-b & b  \notag
\end{vmatrix}\Bigg)=$$
$$=\frac{1}{2}\big(s|t-x|+b|x-t|\big)=\frac{1}{2}(s+b)(x-t)=P(x)-P(t). $$
It follows that $y=s$ and we are done.
\endproof
\section{Limiting shapes and $0$-dimensional subschemes in $\P^2$}
\label{sec:0dim}
Starting from this section we want to use the theory and methods for $0$-dimensional subschemes in $\P^2$ from \cite{CooperHarbourneTeitler} (see for example Remark 2.5.1 for further explanation). We begin with an example which shows how to construct the graph of $\dHF_I$ in some special cases (see \cite[Theorem 4.1.5]{CooperHarbourneTeitler}) for with we can obtain the so-called \textit{reduction vector}.

It is worth mentioning here that in \cite[Section 2.1]{DST15II} the authors described a decomposition algorithm for a finite set of points in the projective plane $\P^2$ and an effective divisor $D$ vanishing to given orders in these points. As a result $D$ can be decomposed into a sum of irreducible curves $C_i$ and a divisor $B(D)$, i.e.,
$$D= \sum_{i=1}^{r} a_i C_i+B(D). $$
This decomposition is called the \textit{Bezout decomposition} and it is uniquely determined (\cite[Theorem 2.6]{DST15II}). The method that we want to use to obtain a reduction vector $u$ is connected to the Bezout decomposition, if we assume that curves the $C_i$ are lines. Each line which appears in the Bezout decomposition will have an impact on the value of the reduction vector, as we will see in the next example.
\begin{example}
	\label{ex:10,8,5,3}
	Consider $4$ distinct lines $L_{10},L_{8},L_{5}$ and $L_{3}$ in $\P^2$ with $10,8,5$ and $3$ general points lying on them respectively. In particular, these points are not intersection points of the lines and non three of them are collinear, unless they lie already on one of the configuration lines. In this example we want to examine the graded sequence of ideals $\ibul=I^{(m)}$, so we assume that in every of distinguished points we have fixed multiplicity $m$.
	
	Assume that number $m$ is least common multiple of numbers $10,8,5$ and $3$, which is $120$. We start by working only with some of the ideals from $\ibul$. We begin with creating the reduction vector $u$ by ``taking'' the line $L_{10}$ as many times, as the sum of all multiplicities  of the points on that line will be equal to the sum of all multiplicities on the line $L_8$. This number, denoted by $k_1$, is equal to
	$$ 10(m-k_1)=8m,\;\; \text{ hence }\;\; k_1=\frac{2}{10}m,$$
	and the reduction vector is
	\[u=
	(\underbrace{10m,10m-10,\ldots,8m}_{\text{$k_1$}},\ldots)
	.\]
	
	In the next step we we need to choose between the lines $L_{10}$ and $L_8$ according to the sum of multiplicities on them (we choose the lines with the higher sum of multiplicities, if they are equal we choose an arbitrary line from this pair of lines).
	Let $k_2$ and $k_3$ be the numbers counting how many times each of lines were chosen until the sum of the multiplicities will be $5m$. As before we compute these numbers using equations
	$$10(\frac{8}{10}m-k_2)=5m, \;\; 8(m-k_3)=5m,$$
	from which $k_2=\frac{3}{10} m$ and $k_3=\frac{3}{8}m. $
	Then we get
	\[u=
	(\underbrace{10m,\ldots,8m}_{\text{$k_1$}},\underbrace{8m-8,\ldots,5m}_{\text{$k_2+k_3$}},\ldots)
	.\]
	We repeat the procedure for lines $L_{10},L_8$ and $L_5$, obtaining the numbers $k_4,k_5$ and $k_6$ from equations
	$$10(\frac{5}{10}m-k_4)=3m,\;\;8(\frac{5}{8}m-k_5)=3m,\;\;5(m-k_6)=3m.$$
	
	\noindent
	To obtain the graph of $\overline{\dHF_{\ibul}}$ we need a few additional observations:
	\begin{itemize}
		\item[1)] It is not hard to see that for $0$-dimension subscheme the values of $\overline{\dHF_{\ibul}}(t)$ are equal to $0$ for $t \geq \widetilde{\reg}(\ibul)$. It is so since for $0$-dimension subscheme the values of $\HF_{I_m}(t)$ are constant if $t \geq \reg(I_m)$.
		\item[2)] Since $(I_m)_t = 0$ for $t < \alpha(I_m)$, we see that $\HF_{I_m}(t) = \binom{t+2}{2}$ for all $t < \alpha(I_m)$ and then
		$$\dHF_{I_m}(t)=\binom{t+2}{2}-\binom{t+1}{2}=t+1.$$
		Thus for $0\leq t\leq \widehat{\alpha}(\ibul)$ we have
		$$ \overline{\dHF_{\ibul}}(t)=\lim_{m \rightarrow \infty} \frac{\Delta\widetilde{HF}_{I_m}(mt)}{m}=t.$$
		\item [3)] The graph of $\dHF_{I_m}$ results from the entries of the reduction vector in the following way. Given the $i$-th entry $a$ in the reduction vector, we mark the points $(i-1,i),(i,i),\ldots,(i+a-2,i)$ on the horizontal line $y=i$. The graph of $\dHF_{I_m}$ is the upper envelope (see \cite{Yates} p.75-80) of the resulting set of marked lattice points.
		\item [4)] If the number $m$ is divisible by at least one of the numbers $10,8,5$ or $3$, the values of all $k_i$ should be rounded up, i.e., $k_1= \lceil \frac{2}{10}m \rceil$, $k_2= \lceil \frac{3}{10}m \rceil$ and etc. It is obvious that $k \leq \lceil k \rceil$, so from the definition of the sets $L_{m,t}(\ibul)$ and therefore from the construction of the set $\Gamma(\ibul,t)$ we see that it is enough to consider only numbers $m$ which are exactly divisible by $10,8,5$ and $3$.
	\end{itemize}
	Given the previous calculations and observations, we obtain
	\[u=
	(\underbrace{10m,\ldots,8m}_{\text{$k_1$}},\underbrace{8m,\ldots,5m}_{\text{$k_2+k_3$}},
	\underbrace{5m,\ldots,3m}_{\text{$k_4+k_5+k_6$}},\ldots,0)
	.\]
	
	Additionally with Corollary \ref{cor:reading gin from dHF graph in P^2}, we conclude that the graph of $\overline{\dHF_{\ibul}}$ and the set $\Gamma(\ibul)$ are as follows
	\begin{figure}[H]
		\begin{subfigure}
		\centering
		\begin{tikzpicture}[line cap=round,line join=round,>=triangle 45,x=1.0cm,y=1.0cm,scale=0.6]
		\clip(-0.4292177589852011,-0.5824101479915452) rectangle (11.224059196617336,5.2695983086680735);
		\draw [line width=0.8pt] (0.,0.)-- (4.,4.);
		\draw [line width=0.8pt] (4.,4.)-- (4.725,1.725);
		\draw [line width=0.8pt] (4.725,1.725)-- (5.875,0.875);
		\draw [line width=0.8pt] (5.875,0.875)-- (8.2,0.2);
		\draw [line width=0.8pt] (8.2,0.2)-- (10.,0.);
		\draw [->,line width=0.8pt] (0.,0.) -- (0.,5.);
		\draw [->,line width=0.8pt] (0.,0.) -- (11.,0.);
		\begin{scriptsize}
		\draw (0.2,4.694545454545452) node[anchor=north west] {$\overline{\dHF_{\ibul}}$};
		\draw (3.088752642706133,4.694545454545452) node[anchor=north west] {$(\widehat{\alpha}(I_{\bullet}),\widehat{\alpha}(I_{\bullet}))$};
		\draw (4.678604651162793,2.377420718816065) node[anchor=north west] {$(3+\sum_{i=1}^6 \frac{k_i}{m},\sum_{i=1}^6 \frac{k_i}{m})$};
		\draw (5.642663847780129,1.6332346723044375) node[anchor=north west] {$(5+\sum_{i=1}^3 \frac{k_i}{m},\sum_{i=1}^3 \frac{k_i}{m})$};
		\draw (7.756828752642709,1.007441860465114) node[anchor=north west] {$(8+\frac{k_1}{m},\frac{k_1}{m})$};
		\draw [fill=black] (0.,0.) circle (0.5pt);
		\draw [fill=black] (4.,4.) circle (1.5pt);
		\draw [fill=black] (4.725,1.725) circle (1.5pt);
		\draw [fill=black] (5.875,0.875) circle (1.5pt);
		\draw [fill=black] (8.2,0.2) circle (1.5pt);
		\draw [fill=black] (10.,0.) circle (1.5pt);		
		\draw (-0.46304439746300236,4.407019027484141) node[anchor=north west] {$4$};
		\draw (-0.46304439746300236,2.407019027484141) node[anchor=north west] {$2$};
		\draw (-0.3446511627906979,0.19560253699788385) node[anchor=north west] {$0$};
		\draw (3.718773784355179,-0.007357293868923736) node[anchor=north west] {$4$};
		\draw (9.718773784355179,-0.007357293868923736) node[anchor=north west] {$10$};
		\draw (7.718773784355179,-0.007357293868923736) node[anchor=north west] {$8$};
		\draw (5.718773784355179,-0.007357293868923736) node[anchor=north west] {$6$};
		\draw (10.718773784355179,-0.007357293868923736) node[anchor=north west] {$t$};
		\end{scriptsize}
		\end{tikzpicture}
		\end{subfigure}
		\begin{subfigure}{}
		\centering
		\begin{tikzpicture}[line cap=round,line join=round,>=triangle 45,x=1.0cm,y=1.0cm,scale=0.6]
		\clip(-0.4292177589852011,-0.5824101479915452) rectangle (11.224059196617336,5.2695983086680735);
		\draw [line width=0.8pt] (0.,0.)-- (4.,4.);
		\draw [line width=0.8pt] (4.,4.)-- (4.725,1.725);
		\draw [line width=0.8pt] (4.725,1.725)-- (5.875,0.875);
		\draw [line width=0.8pt] (5.875,0.875)-- (8.2,0.2);
		\draw [line width=0.8pt] (8.2,0.2)-- (10.,0.);
		\draw [->,line width=0.8pt] (0.,0.) -- (0.,5.);
		\draw [->,line width=0.8pt] (0.,0.) -- (11.,0.);
		\begin{scriptsize}
		\draw (0.2,4.694545454545452) node[anchor=north west] {$\overline{\dHF_{\ibul}}$};
		\draw (3.5792389006342495,4.694545454545452) node[anchor=north west] {$(4,4)$};
		\draw (4.560211416490486,2.7326004228329785) node[anchor=north west] {$\Big(\frac{189}{40},\frac{69}{40}\Big)$};
		\draw (5.761057082452432,1.8192811839323444) node[anchor=north west] {$\Big(\frac{47}{8},\frac {7} {8}\Big)$};
		\draw (7.672262156448204,1.3118816067653256) node[anchor=north west] {$\Big(\frac{41}{5},\frac{1}{5}\Big)$};
		\draw (9.600380549682875,0.7199154334038034) node[anchor=north west] {$(10,0)$};
		\draw (-0.46304439746300236,4.407019027484141) node[anchor=north west] {$4$};
		\draw (-0.46304439746300236,2.407019027484141) node[anchor=north west] {$2$};
		\draw (-0.3446511627906979,0.19560253699788385) node[anchor=north west] {$0$};
		\draw (3.718773784355179,-0.007357293868923736) node[anchor=north west] {$4$};
		\draw (9.718773784355179,-0.007357293868923736) node[anchor=north west] {$10$};
		\draw (7.718773784355179,-0.007357293868923736) node[anchor=north west] {$8$};
		\draw (5.718773784355179,-0.007357293868923736) node[anchor=north west] {$6$};
		\draw (10.718773784355179,-0.007357293868923736) node[anchor=north west] {$t$};
		\draw [fill=black] (4.725,1.725) circle (0.5pt);
		\draw [fill=black] (5.875,0.875) circle (0.5pt);
		\draw [fill=black] (8.2,0.2) circle (0.5pt);
		\end{scriptsize}
		\end{tikzpicture}
		\end{subfigure}
\caption{. The graphs of $\overline{\dHF_{\ibul}}$ for $10,8,5$ and $3$ points on $4$ lines.}
	\end{figure}
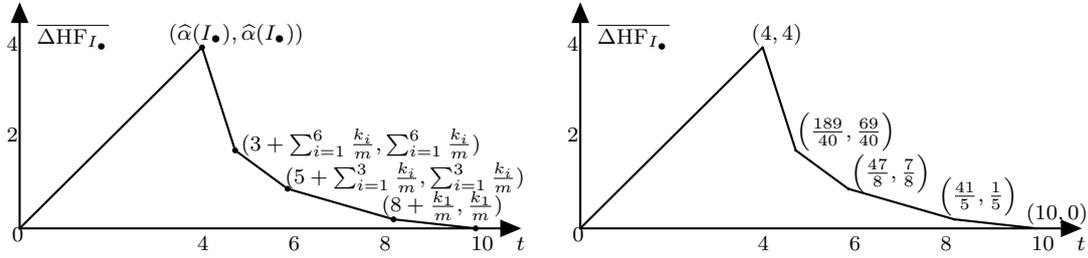

	\begin{figure}[H]
		\centering
		\includegraphics[height=6cm,width=8cm]{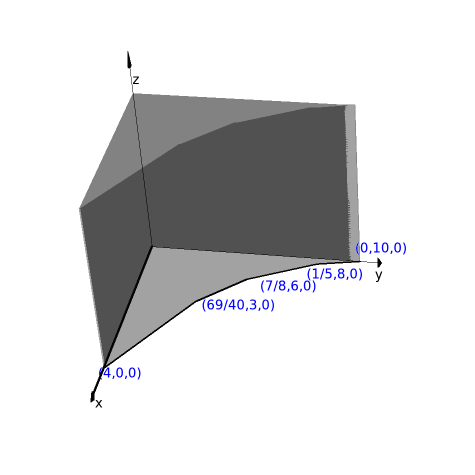}
		\caption{. Shaded area presents the set $\Gamma(\ibul)$.}
	\end{figure}
\end{example}
We can easily see that points chosen in the previous example are very special and this speciality leads to result where the limiting shape consists of some break points. As the next step we are going to describe all sets of points for which we can proceed in the similar way.

\begin{definition}
	\label{def:adapted set}
	Given a set of points $Z \subset \P^2$. We say that a configuration of lines $\mathcal{L}=\{L_1,\ldots,L_n\}$ is \textit{adapted} to the set $Z$, if
	\begin{itemize}
		\item $Z \subseteq L_1\cup\ldots \cup L_n$;
		\item $Z\cap (L_i \cap L_j) = \emptyset $ for all $i<j$;
	\end{itemize}
	Moreover, if $C$ is an irreducible curve of degree $d$ such that $C \cap Z \geq \binom{d+2}{2}$, then $d=1$ and $C \in \mathcal{L}$.
\end{definition}
\begin{remark}
	\label{rmk:adapted sets and blow up}
	If $Z$ is a set of points, such that there exists an adapted configuration of lines for $Z$, then in the blow up $f:\; X \longrightarrow \P^2$
	of $\P^2$ at $Z$, there are only $(-1)$-curves coming from (being proper transforms of) exceptional divisor, lines through pairs of points in $Z$, conics through quintuples of points in $Z$ and possibly more negative curves coming from configuration of lines.
\end{remark}
\begin{remark} There are sets $Z\subset\P^2$ such that there is no adapted configuration to $Z$.
\vspace{0.2cm}	

	\noindent
	\begin{minipage}{0.6\textwidth}
		 For example the set of points presented on Figure \ref{fig: 9dots} is the set of points for which, according to Definition \ref{def:adapted set}, we cannot adapt a configuration of lines. Indeed, all $3$ horizontal and vertical lines must belong to the set $\mathcal{L}$, but this configuration of lines does not fulfil the second condition of Definition \ref{def:adapted set}. Hence it makes sense to state that a set $Z \subset \P^2$ is \textit{adaptable} if there exists a configuration of lines adapted to $Z$.
	\end{minipage}
	\begin{minipage}{0.39\textwidth}
		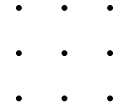
\begin{figure}[H]
			\begin{tikzpicture}[line cap=round,line join=round,>=triangle 45,x=1.0cm,y=1.0cm,scale=0.6]
			\draw [fill=black] (0.,0.) circle (1.5pt);
			\draw [fill=black] (1.,0.) circle (1.5pt);
			\draw [fill=black] (2.,0.) circle (1.5pt);
			\draw [fill=black] (0.,1.) circle (1.5pt);
			\draw [fill=black] (1.,1.) circle (1.5pt);
			\draw [fill=black] (2.,1.) circle (1.5pt);
			\draw [fill=black] (0.,2.) circle (1.5pt);
			\draw [fill=black] (1.,2.) circle (1.5pt);
			\draw [fill=black] (2.,2.) circle (1.5pt);
			\end{tikzpicture}
			\caption{. Not adaptable set of points.}
			\label{fig: 9dots}
		\end{figure}
	\end{minipage}
\end{remark}
\begin{definition}
	\label{def:minimal adapted conf}
	Keeping the notation from Definition \ref{def:adapted set} we say that $\mathcal{L}$ is a \textit{minimal configuration adapted to $Z$} if the number of lines in $\mathcal{L}$ is minimal among all adapted configurations. If on each line $L_i$ we have $a_i$ points of multiplicity $m$ then the \textit{weight of} $L_i$ is the number $a_i \cdot m$.
\end{definition}

\begin{theorem}
	\label{th:n lines dHF}
	Let $Z$ be a adaptable set of points in $\P^2$ and let $\mathcal{L}$ be the minimal adapted configuration to the set $Z$. Denote by $a_i$ number of points on the line $L_i$, for $i\in \{1,\ldots,n\}$. We assume that $a_1> a_2 > \ldots > a_n$. Let $\ibul=\{I(mZ)\}_m$. Then the graph of  $\overline{\dHF_{\ibul}}$ is the chain of segments connecting the points:
	\begin{equation}
	\label{eq:n-gon gin}
	(0,0),\;(n,n),\; \Big( a_i+\sum_{\substack{1\leq j \leq i\\ j\leq k \leq i}} \frac{w_{j,k}}{m}, \sum_{\substack{1\leq j \leq i\\ j\leq k \leq i}} \frac{w_{j,k}}{m} \Big),\; (a_1,0),
	\end{equation}
	with $i\geq2$, where
	\begin{equation}
	\label{eq:w pts n-gon gin}
	\frac{w_{j,k}}{m}:=\frac{a_k-a_{k+1}}{a_j}.
	\end{equation}
\end{theorem}
\proof
This theorem is a generalization of Example \ref{ex:10,8,5,3}. Assume that in every point we have the multiplicity $m$, and that this number is divisible by all numbers $a_i$ for all $i$. In all future consideration, if we will think about any specific line $L$, we sometimes use notation $a(L)$ and $m(L)$ instead of $a_i$ and $m$. Observe that the numbers $w_{j,k}$, in (\ref{eq:w pts n-gon gin}) can be calculated with the following reduction algorithm.

\underline{Reduction:}\begin{itemize}
	\item Step $1$. Pick the line $L_1$ with the weight $a_1 m$. The number $a_1 m$ is the first number of the reduction vector $u$. Then decrease the multiplicities in all points on line $L_1$ by $1$. 
	\item Step $i$.
	Denote by $a(L)$ the number of all points with non-zero multiplicity in points lying on the line $L$, and by $m(L)$ the sum of all multiplicities in these points in every step of reduction. Pick a line $L$ with the highest weight $a(L) m(L)$ (if there is more than one such lines, pick any of them). The number $a(L) m(L)$ is the $i$th number of the reduction vector $u$. Then decrease the multiplicities in all points on chosen line by $1$. 
\end{itemize}

Thus to get the number $w_{1,1}$, we need to solve the equation
$a_1 (m-w_{1,1})=m a_{2},$
from which we get $\frac{w_{1,1}}{m}=\frac{a_1-a_{2}}{a_1}$. We may do it in general for the number $w_{j,k}$, assuming that we already calculated the numbers $w_{j,s}$, for all $s$ such that $j\leq s < k$. Then
$$a_j \big((m-w_{j,j}-w_{j,j+1}-\ldots-w_{j,k-1})-w_{j,k}\big)=m a_{k+1}  ,$$
and from that
$$\frac{w_{j,k}}{m}=1-\frac{w_{j,j}}{m}-\frac{w_{j,j+1}}{m}-\ldots-\frac{w_{j,k-1}}{m}-\frac{a_{k+1}}{a_j}$$
$$=\frac{a_j}{a_j}-\frac{a_j-a_{j+1}}{a_j}-\frac{a_{j+1}-a_{j+2}}{a_j}-\ldots-\frac{a_{k-1}-a_{k}}{a_j}-\frac{a_{k+1}}{a_j}
=\frac{a_k-a_{k+1}}{a_j}. $$
Thus we write the reduction vector as follows
\[u=
(\underbrace{a_1\cdot m,\ldots,a_1 (m-w_{1,1})}_{\text{$w_{1,1}$}},\underbrace{a_2\cdot m,\ldots,a_3\cdot m}_{\text{$w_{1,2}+w_{2,2}$}},\ldots,
\underbrace{a_k\cdot m,\ldots,\ldots,a_{k+1}\cdot m}_{\text{$\sum_{j=1}^k w_{j,k}$}},\ldots,0)
.\]

Using the same argumentation which consists of $4$ general facts, as in Example \ref{ex:10,8,5,3}, we obtain the assertion about points (\ref{eq:n-gon gin}) in $\overline{\dHF_{\ibul}}$.
\endproof

\begin{corollary}
	\label{cor:n break points gin}
	The complement in the first octant of the limiting shape for a graded sequence of ideals $\ibul=\{I(mZ)\}_m$ from Theorem \ref{th:n lines dHF} is a cylinder over the convex hull of points
	$$(0,0),\;(n,0),\; \Big(\sum_{\substack{1\leq j \leq i\\ j\leq k \leq i}} \frac{w_{j,k}}{m},  a_i \Big),\; (0,a_1),$$
	for
	$2 \leq i \leq n.$
\end{corollary}
\proof
This is a simple consequence of Corollary \ref{cor:reading gin from dHF graph in P^2},  definition of $\Gamma(\ibul)$ and property of gins for saturated ideals (see \cite[Proposition 2.21]{Green1}).
\endproof

\begin{theorem}
	\label{th:2lines+1pts dHF}
	Let $L_1$ $L_2$ be two lines. Let $Z$ be the set consisting of the intersection point of these lines and additional $a_i$ points on $L_i$, for $i=1,2$, provided that the condition $a_1 a_2> a_1 + a_2$ is fulfilled. Then the graph of $\overline{\dHF_{\ibul}}$ for $\ibul=\{I(mZ)\}_m$ is the chain of segments connecting the points:
	\begin{itemize}
		\item[i)] $(0,0),\; (2,2),\; \Big(\frac{a_1 a_2+ a_1 + a_2}{a_1+a_2},1 \Big),\; \Big(a_2+1,\frac{a_1-a_2}{a_1}\Big),\; (a_1+1,0),$ for $a_1 > a_2$.
		\item[ii)] $(0,0),\; (2,2),\; \Big(\frac{a_1+2}{2},1 \Big),\; (a_1+1,0),$ if $a_1 = a_2$.
	\end{itemize}
\end{theorem}
\proof
Case i). We reduce the proof of this theorem to the construction of the reduction vector $u$. Assume that in each of the marked points we have the multiplicity $m$ and let us assume that this number is divisible by both numbers $a_1$ and $a_1+a_2$. Since $a_1 > a_2$, the first entries in $u$ are obtained by taking the line $L_1$ so many times, that the sum of all multiplicities on lines $L_1$ and $L_2$ will be equal. Denote this number by $k$. Thus
$(a_1+1)(m-k)=a_2 m+(m-k), $
from which
$k=\frac{a_1-a_2}{a_1}m.$
Now, the sum of all multiplicities on each line is
$$(a_1+1)(m-k)=(a_1+1)\frac{a_2}{a_1}m=(a_2+\frac{a_2}{a_1})m,$$
so
\[u=
\Big(\underbrace{(a_1+1)m,\ldots,(a_2+\frac{a_2}{a_1})m}_{\text{$k$}},\ldots\Big)
.\]

In the next step we need to choose between the lines $L_{1}$ and $L_2$ according to the sum of multiplicities on them (we choose the lines with the higher sum of multiplicities, if they are equal we choose an arbitrary line from this pair of lines). Since $a_1>a_2$ as long as the multiplicity in the point of the intersection of the lines is greater than $0$. This gives $\frac{a_2}{a_1}m$ steps. The multiplicity in each point on line $L_i$ will be reduced by $\frac{a_i}{a_1+a_2}$, so the sum of multiplicity on each line is
$$ a_2(1-\frac{a_2}{a_1+a_2})m = a_1(1-\frac{a_1}{a_1+a_2})m,$$
therefore
\[u=
\Big(\underbrace{(a_1+1)m,\ldots,(a_2+\frac{a_2}{a_1})m}_{\text{$\frac{a_1-a_2}{a_1}m$}},
\underbrace{a_2\cdot m,\ldots,a_2(1-\frac{a_2}{a_1+a_2})m}_{\text{$\frac{a_2}{a_1}m$}},\ldots,0\Big)
.\]
We may turn to the construction of $\overline{\dHF_{\ibul}}$. We have two obvious points, $(a_1+1,0)$ from the first entry in $u$, and $(2,2)$ from $\widehat{\alpha}(\ibul)=2$.
Observe that the sum of numbers
$$\frac{a_1-a_2}{a_1}+(a_2+\frac{a_2}{a_1})=a_2+1,$$
gives $x$-coordinate of another point on the graph of $\overline{\dHF_{\ibul}}$ [see fact 3) in Example \ref{ex:10,8,5,3} for more explanation]. The $y$-coordinate of this point is the number of $k$ first steps scaled by $m$, so we conclude that this point is $\Big(a_2+1,\frac{a_1-a_2}{a_1}\Big)$.

We proceed along the same lines, i.e., if we take the sum of number of steps $ \frac{a_1-a_2}{a_1}m+ \frac{a_2}{a_1}m$, and the number describing the sum of multiplicity after last step $a_2(1-\frac{a_2}{a_1+a_2})m$, we get the number
$$\frac{a_1 a_2+ a_1 + a_2}{a_1+a_2}m.$$
If only this number divided by $m$ is greater than $2=\widehat{\alpha}(\ibul)$, what takes place wherever $a_1 a_2> a_1 + a_2$, then the graph of $\overline{\dHF_{\ibul}}$ consists of one more point $\Big(\frac{a_1 a_2+ a_1 + a_2}{a_1+a_2},1 \Big)$.
The number $1$ is the multiplicity of the point of the intersection of $L_1,L_2$ divided by $m$.
We may check by hand that we do not need any other points to determine $\overline{\dHF_{\ibul}}$ by using the fact that the volume of the body between the graph of $\overline{\dHF_{\ibul}}$ and $t$--axis is equal to $\frac{a_1+a_2+1}{2}$ (see \cite[Lemma 2.15]{Mayes}).
Applying observation 4) from Example \ref{ex:10,8,5,3}, we deduce that we are done.

Case ii). In this case we insert into the previous case $a_2=a_1$ and observe that two points, $\Big(a_2+1,\frac{a_1-a_2}{a_1}\Big)$ and $(a_1+1,0)$, coincide.
\endproof

\begin{corollary}
	\label{cor:2L+1pts gin}
	The complementary of the limiting shape for a graded sequence of ideals $\ibul=\{I(mZ)\}_m$ from Theorem \ref{th:2lines+1pts dHF} is a cylinder over convex hull spanned by points
	\begin{itemize}
		\item[i)] $(0,0),\;(2,0),\;\Big(1,\frac{a_1a_2}{a_1+a_2}\Big),\;\Big(\frac{a_1-a_2}{a_1},a_2+\frac{a_2}{a_1}\Big),\;(0,a_1+1) \text{ for }a_1>a_2,$
		\item[ii)] $(0,0),\;(2,0),\;\Big(1,\frac{a_1a_2}{a_1+a_2}\Big),\;(0,a_1+1) \text{ for }a_1=a_2.$
	\end{itemize}
\end{corollary}
\begin{remark}
	\label{rmk:Mayes_answer}
	Other theorems similar to the previous two can be obtained if we proceed in the same way. We see, for example, from Corollary \ref{cor:n break points gin} that we can construct a limiting shape for a graded sequence of symbolic powers of ideal with any, but finite number of line segments forming its boundary. Conclusions following from these two theorems coincide with an observation made by S. Mayes (see \cite[Observation 5.4]{Mayes2}) and contain a partial answer to Question 5.5 in \cite{Mayes2}.
\end{remark}

\paragraph*{\emph{\textbf{Acknowledgement.}}} I want to warmly thank Tomasz Szemberg and Marcin Dumnicki for introducing me to the problem of general initial ideals and for all inspiring discussions through the whole process of writing this article. I also want to thank Tomasz Szemberg for all the suggestions which helped improve the clarity of this presentation. 

The research of author was partially supported by National Science Centre, Poland, grant 2016/21/N/ST1/01491.

\bigskip
\noindent Grzegorz Malara \\
Instytut Matematyki UP,
Podchor\c a\.zych 2,
PL-30-084 Krak\'ow, Poland
\\
\nopagebreak
\textit{E-mail address:} \texttt{grzegorzmalara@gmail.com}


\addcontentsline{toc}{section}{References}
\begin{thebibliography}{99}
\setstretch{0.8}

\bibitem[BoHa10a]{BoHa1}
Bocci, C., Harbourne, B.: Comparing powers and symbolic powers of ideals,
J. Algebr. Geom. 19 (2010), pp. 399-–417.

\bibitem[BoHa10b]{BoHa2}
Bocci, C., Harbourne, B.:
The resurgence of ideals of points and the containment problem,
Proc. Amer. Math. Soc. 138 (2010), pp. 1175--1190

\bibitem[Cha97]{Chandler}
Chandler, K. A.: Regularity of the powers of an ideal, Communications in Algebra 25(12), (1997) pp. 3773--3776.

\bibitem[Chu81]{Chudnovsky}
Chudnovsky, G. V.: Singular points on complex hypersurfaces and multidimensional Schwarz Lemma,
Seminaire de Theorie des Nombres, Paris 1979-–80, Seminaire Delange-Pisot-Poitou,
Progress in Math vol. 12, M.-J. Bertin, editor, Birkhauser,
Boston-Basel-Stutgart 1981

\bibitem[CHT11]{CooperHarbourneTeitler}
Cooper, S., Harbourne, B., Teitler, Z.: Combinatorial bounds on Hilbert function of fat points in projective space,
J. Pure Appl. Alg. 215 (2011), pp. 2165--2179.

\bibitem[Cut00]{Cutkosky}
Cutkosky, S. D.:
Irrational asymptotic behaviour of Castelnuovo-Mumford regularity,
Journal f\"ur die Reine und Angewandte Mathematik 522 (2000), pp. 93 -- 103.

\bibitem[CuKu11]{CutkoskyKurano}
Cutkosky, S.D., Kurano, K.:
Asymptotic regularity of powers of ideals of points in a weighed projective plane,
Kyoto J. Math. 51 (2011), pp. 25--45.

\bibitem[CHT99]{CutkoskyHerzogTrung}
Cutkosky, S. D., Herzog, J., Trung, N. V.:
Asymptotic behaviour of the Castelnuovo -- Mumford regularity,
Compos. Math. 118 (1999), pp. 243--261.

\bibitem[DGPS]{DGPS}
Decker, W.; Greuel, G.-M.; Pfister, G.; Sch{\"o}nemann, H.: 
\newblock {\sc Singular} {4-1-0} --- {A} computer algebra system for polynomial computations.
\newblock {http://www.singular.uni-kl.de} (2017).

\bibitem[DHST14]{DHST}
Dumnicki, M., Harbourne, B., Szemberg, T., Tutaj-Gasi\'{n}ska, H.:
Linear subspaces, symbolic powers and Nagata type conjectures, Adv. Math. 252 (2014), pp. 471--49

\bibitem[DST16]{DST15II}
Dumnicki, M., Szemberg, T., Tutaj-Gasi\'{n}ska, H.:
Symbolic powers of planar point configurations II.,
J. Pure Appl. Alg. 220 (2016), pp. 2001--2016.

\bibitem[DST15]{DST15}
Dumnicki, M., Szpond, J., Tutaj-Gasi\'nska, H.:
Asymptotic Hilbert Polynomial And Limiting Shape.,
J. Pure Appl. Alg. 219 (2015), pp. 4446--4457.

\bibitem[DSST15]{DSST}
Dumnicki, M., Szemberg, T., Szpond, J.,  Tutaj-Gasi\'{n}ska, H.: Symbolic generic initial systems of star configurations,
Journal of Pure and Applied Algebra 219 (2015), pp. 1073--1081.

\bibitem[Eis95]{Eis}
Eisenbud, D.:
Comutative algebra with a view toward algebraic geometry,
Springer-Verlag, New York, 1995.

\bibitem[Eis05]{Eis2}
Eisenbud, D.:
The Geometry of Syzygies,
Springer-Verlag, New York, 2005.

\bibitem[Fek23]{Fekete}
Fekete, H.:
\"{U}ber die Verteilung der Wurzeln bei gewissen algebraischen Gleichungen mit ganzzahligen Koeffizienten,
Math.Z., 17 (1923), pp. 228-–249.

\bibitem[FKL07]{FKL}
De Fernex, T.; Küronya, A., Lazarsfeld, R.: Higher cohomology of divisors on a projective variety,
Mathematische Annalen, Vol. 337, No. 2, (2007), pp. 443--455.

\bibitem[Gal74]{Gal}
Galligo, A.: A propos du th\'{e}or\`{e}me de pr\'{e}paration de Weierstrass. Fonctions de plusieurs
variables complexes, Lecture Notes in Math., Vol. 409 (1974), Springer, Berlin, pp. 543-–579.

\bibitem[GGP95]{GGP}
Geramita A. V., Gimigliano A., Pitteloud Y.: Graded Betti numbers of some embedded rational
n-folds, Math. Ann. 301, (1995), pp. 363--380. Queen's University Mathematical Preprint
No. 1993-12, October, 1993.

\bibitem[Gre98]{Green1}
Green, M. L.: Generic initial ideals. Six lectures on commutative algebra, Progr. Math., 166, Birkh\"{a}user Verlag, Basel, 1998.,  pp. 119–-186.

\bibitem[GrSt98]{Green2}
Green M., Stillman M.: A tutorial on generic initial ideals. Gr\"{o}bner bases and applications.,
(Linz, 1998), pp. 90--108, London Math. Soc. Lecture Note Ser., 251, Cambridge Univ. Press,
Cambridge, 1998.

\bibitem[HHT02]{HerzogHoaTrung}
Herzog, J., Hoa, L. T., Trung, N. V.: Asymtotic linear bounds for the Castelnuovo-Mumford regularity,
Tras. Amer. Math. Soc. 354 (2002), pp. 1793--1809.

\bibitem[Jer54]{Jerison}
Jerison, M.: A property of extreme points of compact convex sets
Proc. Amer. Math. Soc. 5 (1954), pp. 782--783.

\bibitem[May13]{Mayes2}
Mayes, S.: The Asymptotics of Symbolic Generic Initial Systems of Six Points in $P^2$,
{\tt arXiv:1210.8203}, 2013.

\bibitem[May14]{Mayes}
Mayes, S.: The asymptotic behaviour of symbolic generic initial systems of generic points, J.
Pure Appl. Alg. 218 (2014), pp. 381--390.

\bibitem[ReSi80]{ReedSimon}
Reed, M., Simon, B.:
Methods of Modern Mathematical Physics, vol 1.,
Academic Press, Inc., 1980

\bibitem[Swa97]{Swanson}
Swanson, I.: Powers of Ideals, Primary Decompositions, Artin-Rees Lemma and Regularity,
Math. Ann. 307 (1997), pp. 299--313.

\bibitem[Yat52]{Yates}
Yates, R. C.: Envelopes, A Handbook on Curves and Their Properties,
Ann Arbor, MI: J. W. Edwards, (1952), pp. 75--80. 
\end{thebibliography}
\end{document}